\documentclass[11pt]{article}

\usepackage[a4paper,margin=1in]{geometry}
\usepackage{amsmath,amssymb,amsthm}
\usepackage{booktabs}
\usepackage{tabularx}
\newcolumntype{Y}{>{\raggedright\arraybackslash}X}
\newcolumntype{P}[1]{>{\raggedright\arraybackslash}p{#1}}
\usepackage{multirow}
\usepackage{graphicx}
\usepackage{float}
\usepackage{placeins}
\usepackage{subcaption}
\usepackage{enumitem}
\usepackage{xcolor}
\usepackage[mathlines]{lineno}
\usepackage{hyperref}
\usepackage{tikz}
\usepackage[normalem]{ulem}
\usetikzlibrary{arrows.meta,backgrounds,calc,fit,positioning,shadings,decorations.markings}
\usepackage[scaled=.94]{helvet}

\hypersetup{
colorlinks=true, linkcolor=blue!60!black, citecolor=blue!60!black, urlcolor=blue!60!black
}

\newcommand{\Fourier}{\mathcal{F}}
\newcommand{\methodname}{LithoSAGE}
\newcommand{\methodfullname}{Spatially Adaptive and Geometrically Equivariant Neural Operators for Computational Lithography}
\newtheorem{theorem}{Theorem}

\newif\ifshowlinenumbers
\showlinenumbersfalse

\newcommand{\placeholdergraphic}[2]{%
\IfFileExists{#1}{%
        \includegraphics[width=#2]{#1}%
    }{%
\fbox{%
            \begin{minipage}[c][0.24\textheight][c]{#2}
                \centering
                Figure placeholder\\[0.5em]
                \texttt{\detokenize{#1}}
            \end{minipage}%
        }%
    }%
}

\title{\textbf{\methodname{}: \methodfullname}}
 
\author{%
    Ke Jin\textsuperscript{1},\quad
    Zhenli Xu\textsuperscript{1}\thanks{Email:
    \href{mailto:xuzl@sjtu.edu.cn}{xuzl@sjtu.edu.cn}},\quad 
    Zilin You\textsuperscript{1}\thanks{Email:
    \href{mailto:youzilin@sjtu.edu.cn}{youzilin@sjtu.edu.cn}}
    \\[0.6em]
    \small
    \textsuperscript{1} School of Mathematical Sciences,     Shanghai Jiao Tong University,
    Shanghai 200240, China
}
\date{}

\begin{document}
\maketitle
\ifshowlinenumbers
    \linenumbers
\fi

\begin{abstract}
Operator learning offers an efficient surrogate for lithography simulation, replacing costly numerical solutions with direct model inference. However, existing learning architectures are not designed to simultaneously provide geometric equivariance and adaptability to layout windows with different physical extents, and thus their generalization capability is limited. To address this problem, we propose \methodname, a spatially adaptive and geometrically equivariant neural operator method trained on lithography data. Within its basic Fourier neural operator architecture, spectral weights are expressed as continuous functions of normalized radial frequency through a learnable sum of Gaussians. This representation preserves the geometric properties of the lithography system and incorporates rotation and translation equivariance, spatial interaction decay, and direct evaluation on different layout windows. Experiments predict aerial images from mask layouts and resist images from aerial images, demonstrating the geometric equivariance and spatial adaptivity of \methodname{}. Comparisons with other architectures also assess prediction errors and parameter efficiency.

\end{abstract}

\section{Introduction}
\label{sec:introduction}

Pattern transfer by lithography is a central step in semiconductor manufacturing. Lithography simulation follows the conversion of a mask layout into a printed pattern through an optical process and a photoresist process \cite{mack2008fundamental}. In the optical model, the mask transmission function and illumination source determine the aerial image, the intensity distribution within the photoresist. The resist model then converts this aerial image into a virtual continuous resist image; thresholding or contour extraction gives the printed pattern \cite{mack1985prolith,mack2001lithographic}. Diffraction, partial coherence, illumination conditions, and the projection pupil govern the optical transfer of mask features \cite{hopkins1953diffraction}. The resist model describes the ensuing exposure, post-exposure bake, diffusion, chemical reactions, and development that produce the resist contour \cite{dill1975characterization,mack1995acid}. During calibration and mask optimization, these simulations are repeated for each layout clip and process setting. Accounting for all the physical effects increases their computational cost \cite{lin2018data,yang2022generic}. Practical simulators therefore often use calibrated compact process models, in which empirical nonlinear transformations represent detailed resist effects before thresholding or contour extraction \cite{randall1999variable,fuard2003daim,byers2004lpm,granik2007standard}. The resulting sequence of spatial fields makes lithography simulation a natural application of operator learning.

Learning nonlinear mappings between function spaces from data is the central aim of operator learning \cite{kovachki2024operator,subedi2025operator}. When such a mapping comes from a physical model that is costly to evaluate repeatedly, its learned approximation can provide an efficient surrogate. Representative methods include Graph Kernel Network (GKN) \cite{li2020neural}, the Deep Operator Network (DeepONet) \cite{lu2021learning}, and Fourier neural operator (FNO) \cite{li2021fourier}. GKN approximates kernel integral operators by message passing on graphs. Motivated by the universal approximation theorem for nonlinear operators, DeepONet combines coefficients from a branch network with basis functions from a trunk network. FNO instead parameterizes integral kernels in Fourier space, enabling efficient spectral convolution on regular grids. Later developments use other representations and physical information. The Multiwavelet Transform Operator (MWT) approximates operator kernels with multiwavelets \cite{gupta2021multiwavelet}, while the Wavelet Neural Operator (WNO) learns localized representations through wavelet transforms \cite{tripura2023wavelet}. General Neural Operator Transformer (GNOT) uses heterogeneous normalized attention to handle multiple input functions and irregular meshes \cite{hao2023gnot}. For Helmholtz problems, the Neumann Series Neural Operator (NSNO) uses a U-Net to capture multiscale features \cite{NSNO2024}. The Physics-Informed Neural Operator (PINO) includes governing equations in operator learning \cite{PINO2024}, and PINTO addresses representations on complex geometries \cite{PINTO2025}. Multi-scale Fourier Neural Operator (MscaleFNO) employs parallel scaled FNO branches for highly oscillatory mappings \cite{you2026mscalefno, you2026multiscalefourierneuraloperator}. Other recent methods introduce structured long-range interactions. LSR-Net separates long-range and short-range components and represents the former with an exponential Fourier multiplier \cite{hou2026lsrnet}. E-FNO augments Fourier computation with a nonlocal kernel expressed as a sum of exponentials \cite{liu2027efno}.

The optical imaging mapping and the subsequent resist process mapping place lithography simulation naturally within operator learning \cite{yang2022machine}. Many data-driven lithography methods nevertheless use fixed grids and treat the task as image regression. LithoGAN translates mask patterns into resist patterns with a conditional generative adversarial network (GAN) \cite{ye2019lithogan}. Deep Agile Mask Optimization (DAMO) combines a deep lithography simulator, a deep mask generator, and layout splitting to optimize masks over complete chip layouts \cite{chen2020damo}. To represent nonlocal optical coupling between layout features, more recent models incorporate global interactions. The Dual-Band Optics-Inspired Neural Network (DOINN) combines an optimized Fourier unit that mimics frequency-domain optical computation with a convolutional local path that recovers high-frequency layout details \cite{yang2022generic}. The Convolutional Fourier Neural Operator (CFNO) applies a shared FNO to layout tokens and convolves across tokens for efficient mask optimization at large scale \cite{yang2022largescale}. LithoBench supplies circuit layout tiles and an evaluation framework for benchmarking learning methods in lithography simulation and mask optimization \cite{zheng2023lithobench}.

A practical lithography surrogate needs more than accurate predictions on a fixed layout window. Its response should remain consistent when a pattern changes orientation or position, with mask and source rotated jointly for the optical mapping. It should also operate on windows of different physical extents and reflect the relevant decay behavior. Rotation equivariance means that rotating the input field gives the corresponding rotation of the output, provided the optical or resist process has that rotational symmetry. Translation equivariance means that shifting the input shifts the aerial or resist image by the same amount, subject to boundary effects. Spatial adaptivity is the ability to use the model on layout windows with different physical extents. Decay behavior is a fourth consideration: interactions attenuate with spatial separation in the physical domain, and the spectral response attenuates toward high frequencies or the cutoff in the Fourier domain.

Prior work generally treats some of these requirements without combining all of them. In lithography, symmetry has been used to accelerate aerial image simulation and preserve optical kernel symmetry \cite{yu2008fast,guo2024symmetry}, as well as to build compact process models invariant to geometric transformations \cite{granik2020transformational}. Group-equivariant convolutional neural networks (G-CNNs) and U-Net variants impose rotation or reflection symmetries through group convolutions and shared weights \cite{pmlr-v48-cohenc16,ghyselinck2025effectiveness}. Group-equivariant FNOs apply related constructions to spectral layers that respect rotations, translations, and reflections \cite{pmlr-v202-helwig23a}. With discrete rotation groups, the guarantee covers the finite set of represented rotations rather than arbitrary angles. Other neural operator approaches combine convolution and equivariant attention with FNO architectures \cite{zhao2023enhancing}, employ group averaging \cite{erzmann2024equivariant}, or encode translation, rotation, reflection, and scaling for particular physical problems \cite{wei2026symmetric}. Transfer between computational domains and spatial resolutions has also received attention \cite{ye2024locality,bahmani2025resolution}. These approaches, however, generally do not jointly encode decay behavior and support spatial adaptivity over layout windows with different physical extents. Large-layout lithography strategies typically split a layout into fixed-size blocks, process each block with a model or optimizer, and assemble the outputs \cite{chen2020damo,yang2023scalable}. This improves scalability while retaining a dependence on models or optimizers designed for fixed-size blocks. The present work is motivated by combining rotation and translation equivariance with a continuous spectral representation that supports spatial adaptivity and introduces a structural bias toward the required decay behavior.

We propose \methodname{}, a \textbf{S}patially \textbf{A}daptive and \textbf{G}eometrically \textbf{E}quivariant operator learning method for lithography simulation built on FNO. Spatial adaptivity here means applying the same learned operator to layout windows with different physical extents without retraining. Geometric equivariance concerns rotations and translations. The basic construction expresses FNO spectral weights as continuous functions of normalized radial frequency through a learnable sum-of-Gaussians (SOG) expansion. The expansion takes the place of independently assigned weights at discrete Fourier modes, while Fourier filtering continues to represent nonlocal lithographic interactions. Gaussian sums have been used to approximate interaction kernels for efficient kernel summation, convolution, and spectral computation \cite{liang2020kernel,gao2026fastsog}. SOG-Net represented the Fourier multiplier of each convolution layer for long-range interactions by a sum of Gaussian functions with trainable amplitudes and variances, enabling different long-range behaviors of interatomic potentials \cite{ji2025machine}. In \methodname{}, the Gaussian expansion instead represents spectral weight matrices for lithography operators. Continuity in physical frequency allows the learned function to be evaluated on the Fourier grids associated with different physical windows. The same parameters can therefore be used without defining new spectral weights for each window. Radial dependence shares weights among Fourier modes with equal radial frequency and thereby enforces rotation equivariance. Under periodic boundary conditions, Fourier filtering commutes with translations, preserving translation equivariance. The localized Gaussian components and smooth cutoff taper also regularize the spectral response and promote spatial decay of the resulting nonlocal kernels.

The contributions of this work are as follows:
\begin{itemize}[leftmargin=1.5em]

\item We develop a continuous spectral parameterization using a sum of Gaussian functions. A smooth function of physical frequency represents the spectral response, so the same learned function can be evaluated across layout windows. Localized Gaussian bases and a smooth cutoff taper also encode the required decay behavior.

\item We construct radial spectral weights shared by Fourier modes of equal radial frequency. This sharing yields rotation equivariance and retains the translation equivariance of Fourier multipliers.

\item Experiments verify rotation and translation equivariance for both lithography mappings and demonstrate direct evaluation on layout windows of different physical extents. Comparisons with the other considered architectures assess prediction performance, parameter efficiency, equivariance, and transfer across physical window extents and spatial resolutions. Whereas existing methods typically cover only some of these properties, the proposed model combines them in one formulation.

\end{itemize}

The paper proceeds as follows. Sec.~\ref{sec:physical-background} describes the optical and photoresist processes and formulates the two learning tasks. Sec.~\ref{sec:method} develops \methodname{}, including the spectral parameterization by a sum of Gaussian functions and the associated structural design. Sec.~\ref{sec:experiments} gives the experimental setup and results. Sec.~\ref{sec:conclusion} summarizes the method and the experimental findings.

\section{Lithography Processes and Structural Requirements}
\label{sec:physical-background}
\label{sec:lithography-forward-process}
\label{sec:problem-formulation}

In optical lithography, the mask modulates illumination from a source, and the projection system transfers the resulting field to the wafer. The aerial image is the intensity distribution formed inside the photoresist coating \cite{mack2008fundamental}. This distribution initiates photochemical changes in the resist. Post-exposure bake and development then produce the printed resist pattern.

\begin{figure}[htbp]
    \centering
    \includegraphics[width=0.96\linewidth]{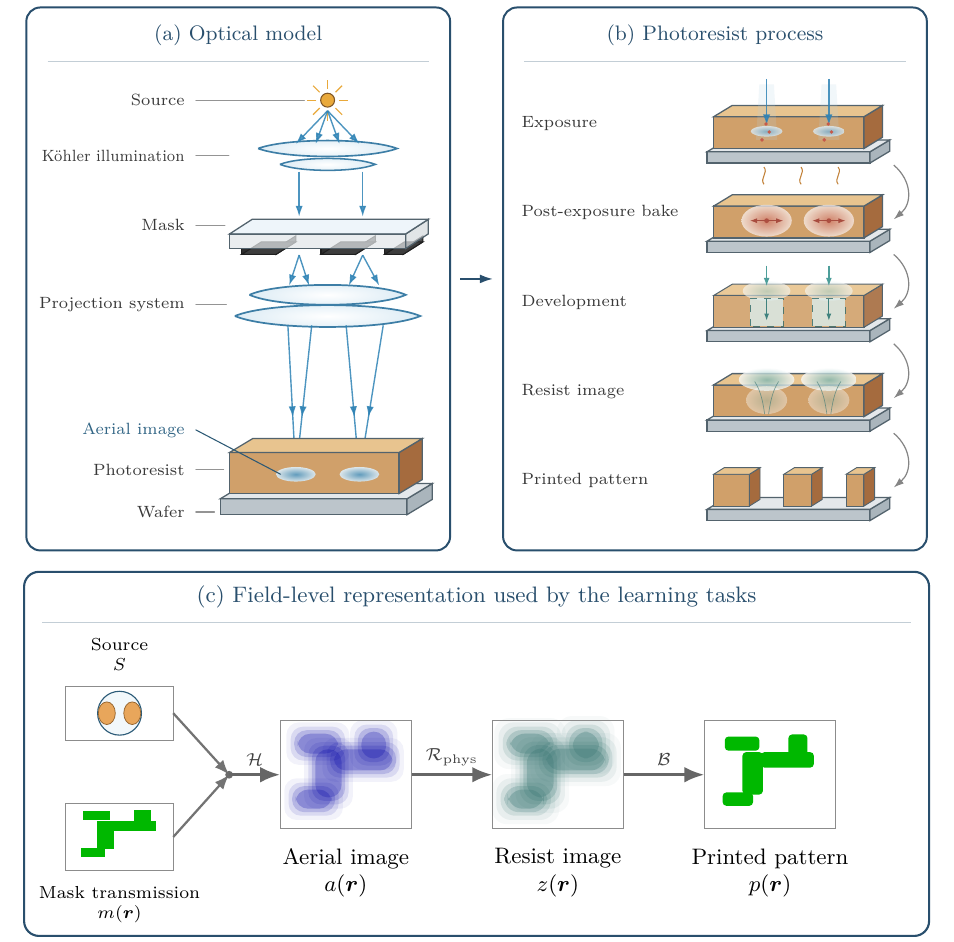}
    \caption{Schematic overview of (a) the optical model, (b) photoresist process, and (c) the field representation used in the learning tasks.}
    \label{fig:lithography-overview}
\end{figure}

The optical imaging sequence is shown in Figure~\ref{fig:lithography-overview}(a). Source light first travels through the Köhler illumination system, which sets the angular distribution incident on the mask. The mask modulates the light according to the layout and generates diffracted light that passes toward the projection system. The finite numerical aperture restricts which diffraction orders the optics can collect, so only part of this light forms the wafer image. Its intensity distribution within the photoresist is the aerial image and sets the spatial exposure for the resist process that follows.

Figure~\ref{fig:lithography-overview}(b) shows the photoresist process following optical imaging. The aerial image drives chemical changes during exposure, with higher intensities producing stronger responses. In a chemically amplified resist (CAR) \cite{mack1995acid,mack2008fundamental}, exposure activates the photoacid generator and generates acid. Diffusion starts during exposure and continues during post-exposure bake (PEB). During PEB, the acid catalyzes reactions that change the local dissolution behavior. Development selectively removes material according to the resulting solubility distribution. In the field representation adopted here, the combined response is described by a continuous \emph{resist image}. Thresholding or contour extraction of that field gives a binary \emph{printed pattern}.

The field representation of this sequence is summarized schematically in Figure~\ref{fig:lithography-overview}(c):
\begin{equation}
    \bigl(m(\boldsymbol{r}),S\bigr)
    \overset{\mathcal{H}}{\longmapsto}
    a(\boldsymbol{r})
    \overset{\mathcal{R}_{\mathrm{phys}}}{\longmapsto}
    z(\boldsymbol{r})
    \overset{\mathcal{B}}{\longmapsto}
    p(\boldsymbol{r}).
    \label{eq:lithography-forward-process}
\end{equation}
The aligned spatial fields share the lateral coordinate $\boldsymbol{r}=(x,y)$ on a computational domain $\Omega\subset\mathbb{R}^2$. We write the mask transmission function, aerial image, resist image, and printed pattern as $m(\boldsymbol{r})$, $a(\boldsymbol{r})$, $z(\boldsymbol{r})$, and $p(\boldsymbol{r})$, respectively. Illumination is specified by the source distribution $S$ in normalized pupil coordinates. The optical mapping, photoresist process, and thresholding or contour extraction are denoted by $\mathcal{H}$, $\mathcal{R}_{\mathrm{phys}}$, and $\mathcal{B}$, respectively. Panel (c) shows the source and mask entering the optical model, followed by the spatial output fields. The physical models for these steps are given below.

\subsection{Optical Model}
\label{sec:optical-imaging}

The optical operator $\mathcal{H}$ in Eq.~\eqref{eq:lithography-forward-process} takes a mask transmission function (MTF) $m(\boldsymbol{r})$ and a continuous source distribution $S(\boldsymbol{s})$ to the aerial image $a(\boldsymbol{r})$ inside the photoresist \cite{goodman2017introduction,mack2008fundamental}. The lateral coordinate is $\boldsymbol{r}=(x,y)\in\Omega$, as defined above, and $m(\boldsymbol{r})$ specifies the mask's spatial transmission. The mask layout describes the arrangement of features; $m(\boldsymbol{r})$ is the transmission field used in the optical calculation. For an ideal binary amplitude mask, $m(\boldsymbol{r})$ is one in transmissive regions and zero in opaque regions. Each coordinate $\boldsymbol{s}=(s_x,s_y)$ identifies an illumination direction in normalized pupil coordinates, with source intensity density $S(\boldsymbol{s})$.

For a fixed direction of illumination $\boldsymbol{s}$, a plane wave represents the field reaching the mask. The scalar thin mask approximation describes its local modulation by the MTF as
\begin{equation}
\begin{aligned}
    U_{\mathrm{in}}(\boldsymbol{r};\boldsymbol{s})
    &=
    e^{-2\pi\mathrm{i}f_c
    \boldsymbol{s}\cdot\boldsymbol{r}}, \\
    U_{\mathrm{mask}}(\boldsymbol{r};\boldsymbol{s})
    &=
    m(\boldsymbol{r})
    U_{\mathrm{in}}(\boldsymbol{r};\boldsymbol{s})
    =
    m(\boldsymbol{r})
    e^{-2\pi\mathrm{i}f_c
    \boldsymbol{s}\cdot\boldsymbol{r}},
\end{aligned}
\label{eq:mask-modulated-field}
\end{equation}
where $\mathrm{i}$ denotes the imaginary unit. The optical cutoff $f_c=\mathrm{NA}/\lambda$ is set by the numerical aperture $\mathrm{NA}$ and wavelength $\lambda$.

Propagation through the projection optics is described by expanding the MTF into spatial-frequency components:
\begin{equation}
m(\boldsymbol{r})
=
\mathcal{F}^{-1}[\widehat{m}](\boldsymbol{r})
=
\int_{\mathbb{R}^2}
\widehat{m}(\boldsymbol{f})
e^{2\pi\mathrm{i}\boldsymbol{f}\cdot\boldsymbol{r}}
\,\mathrm{d}\boldsymbol{f},
\label{eq:mask-fourier-inverse}
\end{equation}
where $\boldsymbol{f}=(f_x,f_y)$ is the spatial frequency coordinate, and $\widehat{m}(\boldsymbol{f})$ is its spectral amplitude. The transform is defined by
\begin{equation}
\widehat{m}(\boldsymbol{f})
=
\mathcal{F}[m](\boldsymbol{f})
=
\int_{\Omega}
m(\boldsymbol{r})
e^{-2\pi\mathrm{i}\boldsymbol{f}\cdot\boldsymbol{r}}
\,\mathrm{d}\boldsymbol{r}.
\label{eq:mask-fourier}
\end{equation}

For illumination along $\boldsymbol{s}$, the oblique incident plane wave shifts the mask spectrum by $f_c\boldsymbol{s}$ relative to the projection pupil. The same dependence on illumination can be expressed through a shifted pupil, $\Pi(\boldsymbol{f}-f_c\boldsymbol{s})$. For an ideal circular pupil, this function is
\begin{equation}
\Pi(\boldsymbol{f})
=
\begin{cases}
1, & \|\boldsymbol{f}\|_2 \leq f_c,\\
0, & \|\boldsymbol{f}\|_2 > f_c.
\end{cases}
\label{eq:pupil}
\end{equation}

The pupil selects the diffraction components that pass through the projection system. Apart from a phase factor of unit magnitude, the coherent field at the image plane for source direction $\boldsymbol{s}$ is
\begin{equation}
E_{\boldsymbol{s}}(\boldsymbol{r};m)
=
\mathcal{F}^{-1}
\left[
\widehat{m}(\boldsymbol{f})
\Pi(\boldsymbol{f}-f_c\boldsymbol{s})
\right](\boldsymbol{r}),
\label{eq:coherent-field}
\end{equation}
where $\mathcal{F}^{-1}$ is the inverse Fourier transform. The intensity contributed by this illumination direction is thus $|E_{\boldsymbol{s}}(\boldsymbol{r};m)|^2$.

The source is normalized by $\int_{\mathbb{R}^2} S(\boldsymbol{s})\,\mathrm{d}\boldsymbol{s}=1$, and $S(\boldsymbol{s})=0$ outside the illumination support. Under K\"ohler illumination, electromagnetic waves from different illumination directions are mutually incoherent. Their intensity contributions can therefore be added, giving the continuous Abbe formulation
\begin{equation}
a(\boldsymbol{r};m,S)
=
\int_{\mathbb{R}^2}
S(\boldsymbol{s})
\left|E_{\boldsymbol{s}}(\boldsymbol{r};m)\right|^2
\,\mathrm{d}\boldsymbol{s}.
\label{eq:abbe-continuous}
\end{equation}
By summing the intensity contributions of individual illumination directions, the Abbe formulation directly represents partially coherent imaging in a computational form. The equivalent Hopkins formulation combines the illumination source and projection pupil in the transmission cross coefficient (TCC), which is widely used in practical lithography simulation \cite{hopkins1953diffraction,mack2008fundamental}. We use the Abbe form because access to individual source points is convenient for optical data generation and for the source-conditioned construction developed later. Let $\mathcal{M}$, $\mathcal{S}$, $\mathcal{A}$, and $\mathcal{Z}$ be the spaces of mask transmission functions, continuous source distributions, aerial images, and resist images, respectively. We use $\theta$ generically for the trainable parameters of each surrogate.

We seek a learned surrogate operator
\begin{equation}
    \widetilde{\mathcal{H}}_{\theta}:
    \mathcal{M}\times\mathcal{S}
    \rightarrow
    \mathcal{A},
    \label{eq:source-conditioned-optical-operator}
\end{equation}
whose predicted aerial image, for an MTF $m(\boldsymbol{r})$ and continuous source distribution $S(\boldsymbol{s})$, satisfies
\begin{equation}
    \widetilde{a}(\boldsymbol{r})
    =
    [\widetilde{\mathcal{H}}_{\theta}(m,S)](\boldsymbol{r})
    \approx
    [\mathcal{H}(m,S)](\boldsymbol{r})
    =
    a(\boldsymbol{r};m,S).
    \label{eq:mask-to-ai-objective}
\end{equation}
The source distribution is specified in normalized pupil coordinates and supplied separately as a conditioning input. It is neither a spatial field on $\Omega$ nor an additional field concatenated with the MTF.

\subsection{Photoresist Process}
\label{sec:photoresist-process}
\label{sec:car-compact-model}

The photoresist process in Figure~\ref{fig:lithography-overview}(b) is described here as a mapping from the aerial image $a(\boldsymbol{r})$ to a continuous resist image $z(\boldsymbol{r})$. In a chemically amplified resist (CAR), light activates the photoacid generator. The resulting acid catalyzes deprotection reactions, altering the local dissolution rate. Acid diffuses during exposure and the subsequent PEB \cite{mack1995acid}. Development then selectively removes material according to its solubility. A detailed simulation can include coupled transport and reaction dynamics and variation through the resist thickness, but requires calibrated material and process parameters \cite{mack2008fundamental}.

For a tractable approximation, we use a two-dimensional CAR compact model consisting of spatial filtering and nonlinear response functions \cite{fuard2003daim,byers2004lpm,randall1999variable,granik2007standard}. The lateral coordinate is $\boldsymbol{r}=(x,y)\in\Omega\subset\mathbb{R}^2$. Instead of resolving the fully coupled resist dynamics, the model treats exposure, acid diffusion, chemical amplification, and development as separately modeled responses in a sequence. An aerial image satisfying $a(\boldsymbol{r})\geq 0$ is mapped to a continuous resist image $z:\Omega\rightarrow[0,1]$ by the following four stages:
\begin{align}
    \text{Exposure response:}\quad
    &H_0(\boldsymbol{r})
    =
    1-\exp\!\left[
        -\alpha_{\mathrm{exp}}D\,a(\boldsymbol{r})
    \right],
    \label{eq:car-exposure}\\
    \text{Acid diffusion:}\quad
    &H_{\mathrm{diff}}(\boldsymbol{r})
    =
    (\kappa_{\sigma}*H_0)(\boldsymbol{r}),
    \label{eq:car-diffusion}\\
    \text{Chemical amplification:}\quad
    &H_{\mathrm{amp}}(\boldsymbol{r})
    =
    1-\exp\!\left[
        -k_{\mathrm{amp}}
        H_{\mathrm{diff}}(\boldsymbol{r})^{n_{\mathrm{amp}}}
    \right],
    \label{eq:car-amplification}\\
    \text{Development response:}\quad
    &z(\boldsymbol{r})
    =
    \Theta_{\beta}\!\left(
        H_{\mathrm{amp}}(\boldsymbol{r})
        -\vartheta_{\mathrm{dev}}
    \right).
    \label{eq:car-resist}
\end{align}
The intermediate fields $H_0$, $H_{\mathrm{diff}}$, and $H_{\mathrm{amp}}$ are the exposure, diffused, and amplified responses, respectively. The normalized exposure dose $D$ scales the aerial image, and $\alpha_{\mathrm{exp}}$ sets the exposure sensitivity. Here $*$ denotes convolution. Convolution with $\kappa_{\sigma}$ models cumulative lateral acid diffusion over exposure and PEB, rather than resolving acid transport in time. A response at one location therefore depends on exposure in neighboring regions as well as the local aerial image intensity. The parameters $k_{\mathrm{amp}}$ and $n_{\mathrm{amp}}$ set the strength and shape of the amplification response; $\vartheta_{\mathrm{dev}}$ is the development threshold. Development uses the smooth threshold $\Theta_{\beta}(u)=1/[1+\exp(-\beta u)]$, with $\beta$ controlling transition sharpness.

We model diffusion with an isotropic Gaussian kernel,
\begin{equation}
    \kappa_{\sigma}(\boldsymbol{r})
    =
    \frac{1}{2\pi\sigma^2}
    \exp\!\left(
        -\frac{\|\boldsymbol{r}\|_2^2}{2\sigma^2}
    \right),
    \label{eq:car-gaussian}
\end{equation}
where $\sigma$ gives the effective diffusion length and $\|\boldsymbol{r}\|_2$ measures Euclidean distance from the kernel center. Dependence of $\kappa_{\sigma}$ only on $\|\boldsymbol{r}\|_2$ gives diffusion no preferred direction. Its decay with physical distance makes nearby aerial image values contribute more to $H_{\mathrm{diff}}$ than distant values. The Fourier transform is also Gaussian, suppressing rapid spatial variations. A threshold of $0.5$ on the resist image defines the printed pattern:
\begin{equation}
    p(\boldsymbol{r})
    =
    \mathbf{1}_{\{z(\boldsymbol{r})\geq 1/2\}}.
    \label{eq:car-printed}
\end{equation}

This compact description does not resolve a full three-dimensional resist simulation. It preserves the spatial coupling from acid diffusion and the local nonlinear responses for exposure, chemical amplification, and development. The diffusion kernel is isotropic, and identical local response functions are used at every position. In this field representation, $m(\boldsymbol{r})$, $a(\boldsymbol{r})$, and $z(\boldsymbol{r})$ follow a periodic boundary convention. The convention serves the numerical treatment of spatial operators; it does not assert physical periodicity of the lithography fields. Isotropic convolution and spatially uniform local responses preserve rotation equivariance. The periodic boundary convention preserves translation equivariance. These features motivate the equivariant architecture presented later.

Denote the CAR compact operator in Sec.~\ref{sec:car-compact-model} by $\mathcal{R}_{\mathrm{CAR}}:\mathcal{A}\rightarrow\mathcal{Z}$. We learn a surrogate $\widetilde{\mathcal{R}}_{\theta}: \mathcal{A}\rightarrow\mathcal{Z}$ whose prediction satisfies
\begin{equation}
    \widetilde{z}(\boldsymbol{r})
    =
    [\widetilde{\mathcal{R}}_{\theta}(a)](\boldsymbol{r})
    \approx
    [\mathcal{R}_{\mathrm{CAR}}(a)](\boldsymbol{r})
    =
    z(\boldsymbol{r}).
    \label{eq:ai-to-ri-objective}
\end{equation}
No explicit source distribution is needed for this mapping, since the aerial image already contains its influence.

\subsection{Structural Requirements}
\label{sec:geometric-requirements}

\begin{figure}[htbp]
    \centering
    \includegraphics[width=0.98\linewidth]{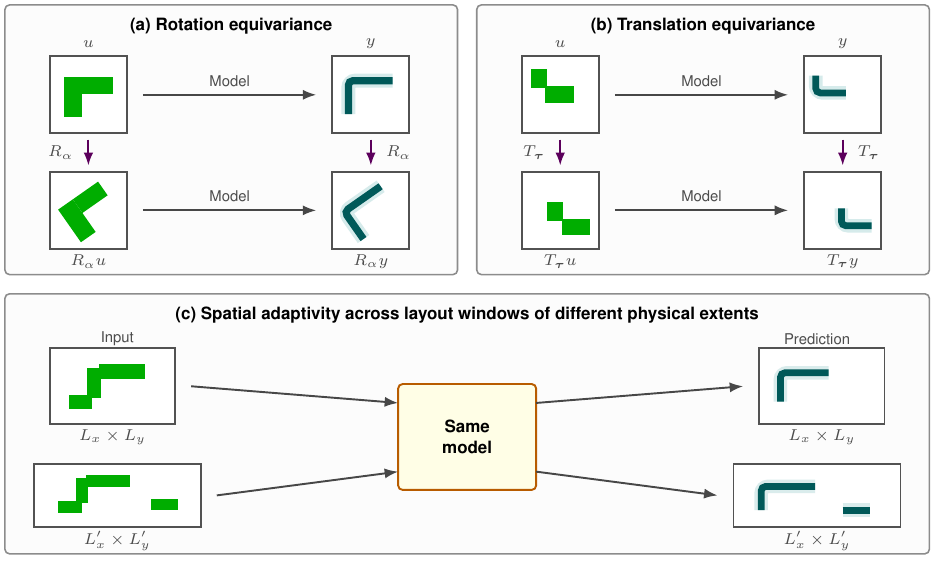}
    \caption{Illustration of the rotation and translation equivariances and spatial adaptivity targeted by the \methodname{} architecture.}
    \label{fig:structural-properties}
    \label{fig:layout-window-adaptivity}
\end{figure}

The physical mapping has geometric, discretization, and interaction structure that a lithography surrogate should preserve. We consider four related properties: rotation equivariance, translation equivariance, spatial adaptivity, and decay behavior. Figure~\ref{fig:structural-properties}(a) depicts the basic rotation relation: rotating the input field should rotate the output in the same way. For the optical mapping, the illumination source rotates together with the mask. Planar rotations constitute the two-dimensional special orthogonal group
\begin{equation}
    SO(2)
    =
    \left\{
    R_\alpha
    =
    \begin{pmatrix}
        \cos\alpha & -\sin\alpha \\
        \sin\alpha &  \cos\alpha
    \end{pmatrix}
    :
    \alpha\in[0,2\pi)
    \right\}.
    \label{eq:so2}
\end{equation}
For a spatial field $u(\boldsymbol{r})$, the rotation acts as
\begin{equation}
    (R_\alpha u)(\boldsymbol{r})
    =
    u(R_\alpha^{-1}\boldsymbol{r}).
    \label{eq:spatial-rotation-action}
\end{equation}
For a continuous source distribution, the analogous action is
\begin{equation}
    (R_\alpha S)(\boldsymbol{s})
    =
    S(R_\alpha^{-1}\boldsymbol{s}).
    \label{eq:source-rotation-action}
\end{equation}
This operation rotates the illumination's angular distribution in normalized pupil coordinates. With an ideal projection system, the continuous Abbe operator obeys the joint equivariance relation
\begin{equation}
    \mathcal{H}(R_\alpha m,R_\alpha S)
    =
    R_\alpha\mathcal{H}(m,S).
    \label{eq:optical-joint-rotation-equivariance}
\end{equation}
Rotation equivariance for mask-to-aerial-image prediction therefore involves both the mask and the illumination source. If an anisotropic dipole source is held fixed while only the mask rotates, Eq.~\eqref{eq:optical-joint-rotation-equivariance} does not generally hold. The continuous optical operator is jointly $SO(2)$ equivariant. In aerial-image-to-resist-image prediction, isotropic diffusion and spatially uniform pointwise response functions have no preferred orientation. A rotation of the aerial image thus produces the same rotation of the resist image.

Figure~\ref{fig:structural-properties}(b) depicts the analogous translation requirement: shifting the input should shift the output by the same amount. For a displacement $\boldsymbol{\tau}=(\tau_x,\tau_y)\in\mathbb{R}^2$, the translation is
\begin{equation}
    T_{\boldsymbol{\tau}}(\boldsymbol{r})
    =
    \boldsymbol{r}+\boldsymbol{\tau}.
    \label{eq:translation}
\end{equation}
Its action on a spatial field $u(\boldsymbol{r})$ is
\begin{equation}
    (T_{\boldsymbol{\tau}}u)(\boldsymbol{r})
    =
    u(\boldsymbol{r}-\boldsymbol{\tau}).
    \label{eq:translation-transformation}
\end{equation}
The inverse argument in Eq.~\eqref{eq:translation-transformation} corresponds to moving the spatial pattern by $\boldsymbol{\tau}$. In optical imaging, source coordinates specify illumination angles in the pupil plane, so they do not move with a spatial translation of the mask. A spatially homogeneous optical model consequently satisfies
\begin{equation}
    \mathcal{H}(T_{\boldsymbol{\tau}}m,S)
    =
    T_{\boldsymbol{\tau}}\mathcal{H}(m,S)
    \label{eq:optical-translation-equivariance}
\end{equation}
provided that the boundary convention preserves translation symmetry. For the resist mapping, diffusion is a function of relative displacement, and the local responses are uniform over the domain. The mapping therefore has translation equivariance under the same boundary convention.

Spatial adaptivity in Figure~\ref{fig:structural-properties}(c) means using a single learned model on layout windows with different physical extents. Lithography simulation involves small windows around individual patterns as well as larger windows that include more surrounding features. The surrogate should return the output on each corresponding window without a separate model for every window size.

Decay behavior, the fourth property, involves attenuation in physical and Fourier spaces. Nonlocal influence should generally decrease with spatial separation. At sufficiently high spatial frequencies, the response should be attenuated and approach zero at the relevant cutoff. Optical imaging combines nonlocal diffraction with a finite passband. Resist diffusion instead combines localized spatial coupling with attenuation of high frequencies. The surrogate should represent these behaviors alongside nonlocal coupling.

\section{Methodology}
\label{sec:method}

\subsection{Fourier Neural Operator Framework}
\label{sec:fourier-framework}

For optical prediction, the input and output are $u(\boldsymbol{r})=m(\boldsymbol{r})$ and $y(\boldsymbol{r})=\widetilde{a}(\boldsymbol{r})$; the illumination source enters separately as a nonspatial conditioning input. For resist prediction, they are $u(\boldsymbol{r})=a(\boldsymbol{r})$ and $y(\boldsymbol{r})=\widetilde{z}(\boldsymbol{r})$.
A standard FNO \cite{li2021fourier} applies a pointwise lifting map, passes the resulting field through Fourier layers, and then applies a pointwise projection:
\begin{equation}
    u(\boldsymbol{r})
    \xrightarrow{P}
    v_0(\boldsymbol{r})
    \xrightarrow{\mathcal{L}_0}
    \cdots
    \xrightarrow{\mathcal{L}_{L-1}}
    v_L(\boldsymbol{r})
    \xrightarrow{Q}
    y(\boldsymbol{r}).
    \label{eq:fno-architecture}
\end{equation}
The map $P$ lifts the scalar input into a multichannel latent field, and $Q$ projects the last latent field to a scalar output. The network contains $L$ Fourier layers: $\mathcal{L}_\ell$ is layer $\ell$, with input latent field $v_\ell$. We use $d$ for the common latent channel width. A Fourier layer takes the form
\begin{equation}
    v_{\ell+1}(\boldsymbol{r})
    =
    \sigma_{\mathrm{act}}\left(
        W_\ell v_\ell(\boldsymbol{r})
        +
        (\mathcal{K}_\ell v_\ell)(\boldsymbol{r})
    \right),
    \label{eq:fno-layer}
\end{equation}
where $W_\ell$ mixes channels through a pointwise linear map, $\sigma_{\mathrm{act}}$ is a pointwise activation, and $\mathcal{K}_\ell$ acts as the Fourier multiplier
\begin{equation}
    (\mathcal{K}_\ell v_\ell)(\boldsymbol{r})
    =
    \Fourier^{-1}\left(
        G_\ell(\boldsymbol{k}) \cdot
        \Fourier(v_\ell)(\boldsymbol{k})
    \right)(\boldsymbol{r}).
    \label{eq:fno-fourier-operator}
\end{equation}
Here $\Fourier$ and $\Fourier^{-1}$ are the discrete Fourier transform and its inverse, respectively. The discrete frequency index is $\boldsymbol{k}=(k_x,k_y)$, and $G_\ell(\boldsymbol{k})$ is the spectral weight matrix for layer $\ell$.

Write the latent-field Fourier coefficient as $\widehat{v}_\ell(\boldsymbol{k}) = \Fourier(v_\ell)(\boldsymbol{k})$ at discrete index $\boldsymbol{k}=(k_x,k_y)$. The retained indices form a set $\Lambda$. In layer $\ell$, let $C_{\mathrm{in}}$ and $C_{\mathrm{out}}$ be the input and output channel counts, respectively. For every $\boldsymbol{k}\in\Lambda$, a standard FNO learns a spectral matrix $G_\ell(\boldsymbol{k})\in\mathbb{C}^{C_{\mathrm{out}}\times C_{\mathrm{in}}}$. The output Fourier coefficient in channel $o$ is computed by
\begin{equation}
    \widehat{(\mathcal{K}_\ell v_\ell)}_o(\boldsymbol{k})
    =
    \sum_{i=1}^{C_{\mathrm{in}}}
    G_{\ell,o,i}(\boldsymbol{k})
    \widehat{v}_{\ell,i}(\boldsymbol{k}),
    \qquad
    \boldsymbol{k}\in\Lambda.
    \label{eq:standard-fno}
\end{equation}
The indices $i$ and $o$ refer to input and output channels. At frequency index $\boldsymbol{k}$, the entry $G_{\ell,o,i}(\boldsymbol{k})$ transfers the Fourier coefficient from input channel $i$ to output channel $o$. Output coefficients outside the retained set $\Lambda$ are zero.

Fourier multiplication is convolution in physical space, allowing each output location to couple with information throughout the spatial domain. For a rectangular window of side lengths $L_x$ and $L_y$, the discrete Fourier mode index $\boldsymbol{k}=(k_x,k_y)\in\mathbb{Z}^2$ corresponds to angular spatial frequency
\begin{equation}
    (\xi_x,\xi_y)
    =
    \left(
        \frac{2\pi k_x}{L_x},
        \frac{2\pi k_y}{L_y}
    \right).
    \label{eq:physical-frequency}
\end{equation}
\subsection{Continuous Spectral Parameterization}
\label{sec:continuous-spectral-parameterization}

The physical extent of a window determines which physical frequencies its Fourier grid samples. Standard FNO assigns independent spectral weights to a fixed collection of discrete modes. These weights do not enforce equal responses at frequencies related by rotation. If the window extent changes, retaining the same indices changes the physical frequency range of the layer. Keeping that range fixed would instead require a different number of modes, including indices for which no weights have been learned. The construction below uses continuous weights in physical frequency to accommodate different windows and radial sharing to enforce rotation equivariance.

\paragraph{Physical passband.}
Standard FNO fixes the number of retained discrete modes to choose its spectral range. In \methodname{}, a physical-frequency cutoff $\xi_c$ is specified instead, and the retained indices for each window are
\begin{equation}
    \Lambda(L_x,L_y)
    =
    \left\{
    (k_x,k_y)\in\mathbb{Z}^2:
    \sqrt{
    \left(\frac{2\pi k_x}{L_x}\right)^2
    +
    \left(\frac{2\pi k_y}{L_y}\right)^2
    }
    \leq \xi_c
    \right\}.
    \label{eq:physical-frequency-index-set}
\end{equation}
The cutoff $\xi_c$ thus defines a common physical passband, while the number of retained modes varies with $L_x$ and $L_y$. A larger window samples physical frequency more finely and contains more discrete modes in that passband. Its retained mode set adapts to the window extent, with the cutoff and learned spectral parameters unchanged.

\paragraph{Rotation equivariance.}
For the basic Fourier layer, rotation equivariance constrains spectral weights to a radial form. We use angular spatial frequency $\boldsymbol{\xi}=(\xi_x,\xi_y)$ in this analysis; it is related to the optical spatial frequency $\boldsymbol{f}$ by $\boldsymbol{\xi}=2\pi\boldsymbol{f}$. Consider a continuous Fourier layer from input $u\in L^2(\mathbb{R}^2;\mathbb{C}^{C_{\mathrm{in}}})$ to output $y\in L^2(\mathbb{R}^2;\mathbb{C}^{C_{\mathrm{out}}})$, with a bounded, continuous spectral weight matrix $G:\mathbb{R}^2\rightarrow \mathbb{C}^{C_{\mathrm{out}}\times C_{\mathrm{in}}}$. Rotations are assumed to affect spatial coordinates only, leaving the channel components unchanged.

\begin{theorem}[Rotation equivariance of a Fourier layer]
\label{thm:radial-equivariance}
Suppose that
\begin{equation}
\widehat{y}(\xi_x,\xi_y)
=
G(\xi_x,\xi_y)\widehat{u}(\xi_x,\xi_y).
\label{eq:continuous-fourier-layer}
\end{equation}
The layer is equivariant under all planar rotations if and only if the spectral weight $G(\xi_x,\xi_y)$ depends only on the radial frequency $\sqrt{\xi_x^2+\xi_y^2}$. Equivalently, the spectral weight is constant on every circle centered at the origin.
\end{theorem}

Theorem~\ref{thm:radial-equivariance} characterizes $SO(2)$ equivariance of the continuous Fourier layer. It supplies the radial constraint in \methodname{}: equal radial frequencies share a spectral weight. The result holds for real as well as complex weights. If the operator must map real-valued fields to real-valued fields, coefficients at opposite frequencies must also obey conjugate symmetry. Radiality assigns equal weights to these opposite frequencies, so satisfying that symmetry requires real spectral weights. Appendix~\ref{app:radial-proof} proves the theorem, and Appendix~\ref{app:real-spectral-weights} gives the conjugate symmetry conditions for the two mappings.

In the continuous construction, lifting, projection, pointwise channel mixing, and activation use identical maps at every spatial position, preserving rotation and translation equivariance. For the resist mapping, radial weights in each Fourier layer therefore give rotation equivariance of the full network. For the optical mapping, the radial response is measured relative to source coordinates, as defined below; the resulting relation is equivariance under joint rotation of the mask and source. No analogous constraint on spectral weights is needed for translation, because Fourier multipliers already commute with translations.

\paragraph{Gaussian spectral weights.}
Define the radial frequency as $\rho=\sqrt{\xi_x^2+\xi_y^2}$ and its normalized coordinate over the common physical passband as $\nu(\xi_x,\xi_y) = \rho/\xi_c$. The spectral weight for layer $\ell$, output channel $o$, and input channel $i$ is parameterized by
\begin{equation}
    G_{\ell,o,i}(\xi_x,\xi_y)
    =
    \sum_{b=1}^{N_b}
    \gamma_{\ell,o,i,b}
    \exp\!\left(
        -\frac{\left(\nu(\xi_x,\xi_y)-\mu_b\right)^2}
        {\delta_{\ell,b}^2}
    \right).
    \label{eq:radial-basis}
\end{equation}
The expansion uses $N_b$ Gaussian basis functions. Basis $b$ has center $\mu_b$ and positive width $\delta_{\ell,b}>0$ in layer $\ell$. The coefficient $\gamma_{\ell,o,i,b}$ weights this basis from input channel $i$ to output channel $o$ in layer $\ell$. Increasing $N_b$ gives a more flexible representation of the radial spectral response. Centers $\mu_b$ and widths $\delta_{\ell,b}$ may be fixed or trained, whereas coefficients are learned for every layer and channel pair. The coordinate $\nu(\xi_x,\xi_y)$ depends on $(\xi_x,\xi_y)$ only through $\rho(\xi_x,\xi_y)$. Hence Eq.~\eqref{eq:radial-basis} is radial and satisfies Theorem~\ref{thm:radial-equivariance}.

For the optical mapping from mask and illumination source to aerial image, we extend the radial construction using frequencies relative to individual source points. The source-conditioned form represents illumination by a weighted set $S_Q=\{(\boldsymbol{s}_q,w_q)\}_{q=1}^{N_s}$. Here $N_s$ is the number of source points, $\boldsymbol{s}_q=(s_{q,x},s_{q,y})$ gives the normalized pupil coordinate of point $q$, and $w_q$ is a nonnegative weight. The weights are normalized by $\sum_{q=1}^{N_s}w_q=1$. The normalized radial coordinate for source point $q$ is
\begin{equation}
    \nu_q(\xi_x,\xi_y)
    =
    \sqrt{
        \left(
            \frac{\xi_x}{\xi_c}
            -
            \alpha_s s_{q,x}
        \right)^2
        +
        \left(
            \frac{\xi_y}{\xi_c}
            -
            \alpha_s s_{q,y}
        \right)^2
    },
    \label{eq:source-conditioned-radial-coordinate}
\end{equation}
where the dimensionless factor $\alpha_s$ maps normalized pupil coordinates to normalized angular spatial frequency coordinates.

Conditioning the spectral weight on $S_Q$ gives
\begin{equation}
    G_{\ell,o,i}(\xi_x,\xi_y;S_Q)
    =
    \sum_{q=1}^{N_s}
    w_q
    \sum_{b=1}^{N_b}
    \gamma_{\ell,o,i,b}
    \exp\!\left(
        -\frac{
        \left(
            \nu_q(\xi_x,\xi_y)-\mu_b
        \right)^2
        }{\delta_{\ell,b}^2}
    \right).
    \label{eq:source-conditioned-radial-basis}
\end{equation}
The indices $\ell$, $o$, and $i$ and the parameters $N_b$, $\gamma_{\ell,o,i,b}$, $\mu_b$, and $\delta_{\ell,b}$ retain their meanings from Eq.~\eqref{eq:radial-basis}. A basis function is radial in the frequency displacement from one source point. Summing these functions with source weights can nevertheless produce an anisotropic spectral response determined by the illumination. Jointly rotating the mask and source coordinates leaves the relative frequency distances unchanged, so the Fourier layer is equivariant under joint rotation of the mask and source. For the aerial-image-to-resist-image mapping, no explicit source enters and Eq.~\eqref{eq:radial-basis} is used. With the window-dependent mode set $\Lambda(L_x,L_y)$ in Eq.~\eqref{eq:physical-frequency-index-set}, either spectral form can be evaluated on the Fourier grids of square and rectangular windows with different physical extents using the same learned parameters. For the form without an explicit source, radiality together with conjugate symmetry requires real spectral weights.

\begin{figure}[t]
    \centering
    \includegraphics[width=0.93\linewidth]{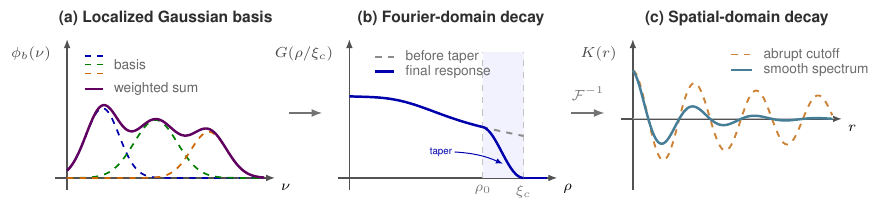}
    \caption{Schematic illustration of the localized Gaussian basis, smooth spectral cutoff taper, and associated spatial kernel decay.}
    \label{fig:sog-decay-mechanism}
\end{figure}

A finite sum of Gaussian radial basis functions can approximate a continuous function on a compact interval arbitrarily well, provided that sufficiently many components and suitable centers and widths are available \cite{park1991universal}. Accordingly, the parameterization in Eq.~\eqref{eq:radial-basis} can, in principle, approximate continuous radial spectral responses over the retained frequency interval as $N_b$ increases and suitable values of $\mu_b$ and $\delta_b$ are provided. For the spectral matrices in \methodname{}, the approximation applies entry by entry. Omitting the layer index for clarity, write the $b$-th Gaussian basis as
\begin{equation}
\phi_b(\nu)
=
\exp\!\left(
-\frac{(\nu-\mu_b)^2}{\delta_b^2}
\right).
\label{eq:gaussian-radial-basis}
\end{equation}
The bases are smooth and localized near $\mu_b$. Their learned weighted sum gives a smooth, continuous spectral response, as in Figure~\ref{fig:sog-decay-mechanism}(a). Near the physical cutoff we apply a cosine taper. It starts at normalized radial frequency $\nu_0$, corresponding to physical radial frequency $\rho_0=\nu_0\xi_c$, and falls smoothly to zero over $\rho_0<\rho<\xi_c$; see Figure~\ref{fig:sog-decay-mechanism}(b). Denote the associated radial spatial kernel by $K_{\mathrm{sp}}(r_{\mathrm{sep}})$, with scalar separation $r_{\mathrm{sep}}$. Abrupt truncation at $\xi_c$ introduces a spectral-edge discontinuity. The smooth taper suppresses it and reduces oscillatory sidelobes in $K_{\mathrm{sp}}(r_{\mathrm{sep}})$, promoting more localized interactions \cite{harris1978windows}, as sketched in Figure~\ref{fig:sog-decay-mechanism}(c). These curves are schematic. The learned response need not be globally monotonic, and the parameterization imposes no particular spatial decay law.

\begin{figure}[htbp]
    \centering
    \includegraphics[width=0.96\linewidth]{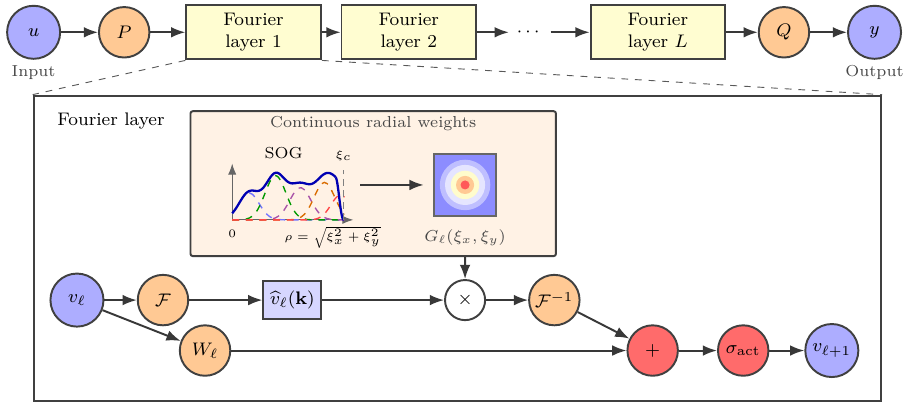}
    \caption{Overall architecture of \methodname{} and the spectral parameterization and computation in one Fourier layer.}
    \label{fig:sogfno-architecture}
\end{figure}

Figure~\ref{fig:sogfno-architecture} summarizes the basic \methodname{} architecture and the operations in one Fourier layer. The lifting map $P$ converts the input field to a latent field, a stack of $L$ Fourier layers processes it, and $Q$ projects it to the output. The upper path shows how continuous radial spectral weights are represented by a sum of Gaussian functions. The lower path applies these weights to Fourier coefficients of $v_\ell$. The Fourier and pointwise outputs are added, and $\sigma_{\mathrm{act}}$ produces the next field $v_{\ell+1}$. For optical prediction, the illumination source additionally conditions the spectral weights, as in Eq.~\eqref{eq:source-conditioned-radial-basis}.

\section{Numerical Results}
\label{sec:experiments}

We compare \methodname{} against FNO \cite{li2021fourier}, DOINN \cite{yang2022generic}, and CFNO \cite{yang2022largescale}. The comparisons cover reference-window prediction errors and parameter counts, geometric equivariance, and transfer across physical window extents and spatial resolutions. The comparisons follow the common data splits and training protocol in Sec.~\ref{sec:experimental-setup}. DOINN and CFNO need input dimensions divisible by their pooling or patch factors. Where necessary, we apply periodic padding on the evaluation grid and crop the output to the original extent. Inputs and targets are never resized.

\subsection{Experimental Protocol}
\label{sec:experimental-setup}

The two tasks in Sec.~\ref{sec:problem-formulation} are mask-to-aerial-image (mask-to-AI) prediction and aerial-image-to-resist-image (AI-to-RI) prediction. We evaluate \methodname{} on both. The first approximates the optical operator of Sec.~\ref{sec:optical-imaging}; the second approximates the CAR compact resist operator of Sec.~\ref{sec:car-compact-model}. The experiments examine the structural properties targeted by the architecture. 

\paragraph{Data generation.} We construct ViaSet and MetalSet from the via and metal layer subsets, respectively, of LithoBench \cite{zheng2023lithobench}. Each dataset uses 1,100 binary mask layouts, divided into 800 training, 100 validation, and 200 test samples with a fixed split. An original binary layout is denoted by $\chi(\boldsymbol{r})\in\{0,1\}$ on the reference square window of side length $L_{\mathrm{ref}}$. Before optical simulation, we low-pass filter $\chi(\boldsymbol{r})$ by retaining frequencies with $\|\boldsymbol{f}\|_2\leq f_{\mathrm{LP}}$. This gives the continuous mask transmission field $m(\boldsymbol{r})$. Area averaging then resamples $m(\boldsymbol{r})$ onto the task grid; small overshoots beyond $[0,1]$ may remain. We simulate the aerial-image target $a(\boldsymbol{r})$ directly on the mask-to-AI grid. The CAR model and thresholding subsequently give the resist image $z(\boldsymbol{r})$ and printed pattern $p(\boldsymbol{r})$. Appendix~\ref{app:physical-data-configuration} lists $L_{\mathrm{ref}}$, $f_{\mathrm{LP}}$, the task grids, and the physical simulation parameters.

ViaSet uses layouts from the first via layers of two benchmark digital circuits. A horizontal or vertical dipole source is assigned to each mask and represented by the quadrature $S_Q$ with $N_s$ points. The two orientations occur equally often in every training, validation, and test split. The optical simulator in TorchResist \cite{wang2025torchresist} generates aerial images with the wavelength $\lambda$, numerical aperture $\mathrm{NA}$, and pupil model in Appendix~\ref{app:physical-data-configuration}. Resist targets use fixed CAR material parameters and the normalized exposure dose $D$ specified there for ViaSet.

MetalSet uses synthesized metal layer layouts with the same split and spatial discretization as ViaSet. Its aerial images use the same source representation and optical parameters $(\lambda,\mathrm{NA})$. With the optical settings held fixed, layout geometry is the main difference between the mask-to-AI datasets. ViaSet has sparse via patterns, while MetalSet has denser patterns with longer, more connected metal features. Resist generation uses the same CAR material parameters, but $D$ takes the MetalSet value in Appendix~\ref{app:physical-data-configuration}. The datasets thus describe the same resist material at different exposure conditions and define different AI-to-RI mappings.

\begin{figure}[htbp]
    \centering
    \includegraphics[width=0.90\linewidth]{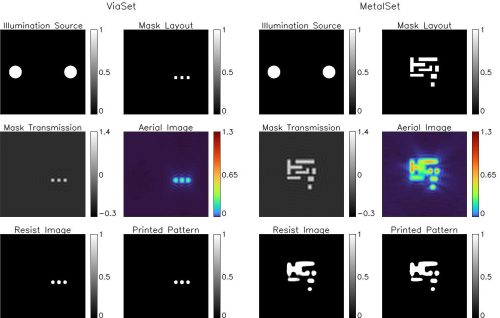}
    \captionsetup{width=0.90\linewidth}
    \caption{Representative generated samples for ViaSet and MetalSet. For each dataset, the panels show the illumination source, mask layout, mask transmission, aerial image, resist image, and printed pattern.}
    \label{fig:dataset-samples}
\end{figure}

Representative ViaSet and MetalSet samples appear in Figure~\ref{fig:dataset-samples}. Mask Layout labels the original binary field $\chi(\boldsymbol{r})$; Mask Transmission labels the continuous field $m(\boldsymbol{r})$ after low-pass filtering. The simulated fields reflect the different geometries. Aerial and resist responses are more localized in ViaSet and extend more broadly along connected features in MetalSet. The resist images additionally reflect the different exposure conditions.

\FloatBarrier

\paragraph{Training protocol.} Equivariance tests and architecture comparisons use models trained on the reference window with the common splits and rotation augmentation. Every model in the architecture comparison has the same training window and epoch budget. Checkpoints are chosen using only the validation split. The spatial adaptivity study uses separate models trained jointly on the reference window and a larger square window, as specified in Sec.~\ref{sec:spatial-adaptivity-results}. Model architectures and training settings are listed in Appendices~\ref{app:lithosage-configuration} and \ref{app:training-configuration}, respectively.

\paragraph{Metrics.} For a reference field and its prediction, write $y,\widetilde{y}\in\mathbb{R}^{N}$, where $N$ counts pixels. We measure continuous aerial and resist image accuracy with the root mean square error (RMSE) and relative $L_2$ error:
\begin{equation}
    \mathrm{RMSE}(\widetilde{y},y)
    =
    \left(
    \frac{1}{N}\sum_{n=1}^{N}
    |\widetilde{y}_n-y_n|^2
    \right)^{1/2},
    \qquad
    \mathrm{Rel.}\ L_2(\widetilde{y},y)
    =
    \frac{\|\widetilde{y}-y\|_2}{\|y\|_2}.
    \label{eq:prediction-errors}
\end{equation}
The maximum absolute error, $\|\widetilde{y}-y\|_{\infty}=\max_n|\widetilde{y}_n-y_n|$, measures the largest pointwise difference in an output field. For AI-to-RI prediction, both the reference and predicted resist images are thresholded at $0.5$. We evaluate the printed patterns with intersection over union (IOU) and PA. Denoting the predicted and reference printed pixel sets by $\mathcal{I}_{\mathrm{pred}}$ and $\mathcal{I}_{\mathrm{ref}}$, respectively, gives
\begin{equation}
    \mathrm{IOU}
    =
    \frac{
        |\mathcal{I}_{\mathrm{pred}}\cap\mathcal{I}_{\mathrm{ref}}|
    }{
        |\mathcal{I}_{\mathrm{pred}}\cup\mathcal{I}_{\mathrm{ref}}|
    },
    \qquad
    \mathrm{PA}
    =
    \frac{
        |\mathcal{I}_{\mathrm{pred}}\cap\mathcal{I}_{\mathrm{ref}}|
    }{
        |\mathcal{I}_{\mathrm{ref}}|
    }.
    \label{eq:iou-pa}
\end{equation}
IOU accounts for extra printed regions as well as missing ones. PA instead measures what fraction of reference printed pixels the prediction recovers. Both metrics range over $[0,1]$, with higher values indicating closer agreement.

For equivariance under a spatial transformation $g$, we compare two outputs: the prediction from the transformed input, $F_\theta(gx)$, and the transformed prediction from the original input, $gF_\theta(x)$. Here $F_\theta$ is the trained model. We quantify the difference with RMSE, relative $L_2$ error, and maximum absolute error. For AI-to-RI prediction, we also threshold both resist images at $0.5$ and compute the IOU of the resulting printed patterns.

\subsection{Prediction Accuracy and Parameter Efficiency}
\label{sec:prediction-comparison}

\label{sec:supplementary-comparison}

\begin{table}[H]
    \centering
    \caption{Prediction errors and parameter counts on the $2048\,\mathrm{nm}\times2048\,\mathrm{nm}$ reference window.}
    \label{tab:supplementary-model-comparison}
    \scriptsize
    \setlength{\tabcolsep}{2pt}
    \begin{tabular}{lllccccc}
        \toprule
        Task & Dataset & Model & Param. ($10^6$) & Output RMSE $\downarrow$ & Output Rel. $L_2$ $\downarrow$ & IOU $\uparrow$ & PA $\uparrow$ \\
        \midrule
        \multirow{8}{*}{Mask$\to$AI}
        & \multirow{4}{*}{ViaSet} & FNO & 4.745 & $7.569{\times}10^{-4}$ & $0.0537$ & -- & -- \\
        & & DOINN & 0.301 & $1.640{\times}10^{-3}$ & $0.1163$ & -- & -- \\
        & & CFNO & 1.577 & $1.445{\times}10^{-3}$ & $0.1024$ & -- & -- \\
        & & \methodname{} & \textbf{0.089} & $4.774{\times}10^{-4}$ & $0.0338$ & -- & -- \\
        \cmidrule(lr){2-8}
        & \multirow{4}{*}{MetalSet} & FNO & 4.745 & $2.362{\times}10^{-3}$ & $0.0293$ & -- & -- \\
        & & DOINN & 0.301 & $5.413{\times}10^{-3}$ & $0.0671$ & -- & -- \\
        & & CFNO & 1.577 & $7.818{\times}10^{-3}$ & $0.0970$ & -- & -- \\
        & & \methodname{} & \textbf{0.089} & $1.998{\times}10^{-3}$ & $0.0248$ & -- & -- \\
        \midrule
        \multirow{8}{*}{AI$\to$RI}
        & \multirow{4}{*}{ViaSet} & FNO & 4.740 & $3.780{\times}10^{-4}$ & $0.0065$ & $0.9966$ & $0.9984$ \\
        & & DOINN & 0.291 & $1.431{\times}10^{-3}$ & $0.0245$ & $0.9884$ & $0.9921$ \\
        & & CFNO & 1.568 & $1.762{\times}10^{-3}$ & $0.0301$ & $0.9854$ & $0.9910$ \\
        & & \methodname{} & \textbf{0.089} & $7.640{\times}10^{-4}$ & $0.0131$ & $0.9927$ & $0.9949$ \\
        \cmidrule(lr){2-8}
        & \multirow{4}{*}{MetalSet} & FNO & 4.740 & $9.060{\times}10^{-4}$ & $0.0055$ & $0.9978$ & $0.9986$ \\
        & & DOINN & 0.291 & $4.139{\times}10^{-3}$ & $0.0251$ & $0.9917$ & $0.9951$ \\
        & & CFNO & 1.568 & $4.826{\times}10^{-3}$ & $0.0292$ & $0.9866$ & $0.9978$ \\
        & & \methodname{} & \textbf{0.089} & $1.457{\times}10^{-3}$ & $0.0088$ & $0.9964$ & $0.9980$ \\
        \bottomrule
    \end{tabular}
\end{table}

Reference-window prediction errors and parameter counts are reported in Table~\ref{tab:supplementary-model-comparison}. With $0.089$ million parameters, \methodname{} uses about $1.9\%$ of FNO's parameter count. Its relative $L_2$ errors are $0.0338$ and $0.0248$ for mask-to-AI prediction on ViaSet and MetalSet, respectively, and $0.0131$ and $0.0088$ for AI-to-RI prediction. Printed-pattern IOU and PA both exceed $0.99$ on the two datasets. These results demonstrate accurate aerial-image and resist-image prediction and high printed-pattern agreement with a compact parameterization.

\subsection{Rotation and Translation Equivariance}
\label{sec:equivariance-results}

We assess rotation and periodic translation equivariance of \methodname{} on both tasks. Each test compares two outputs: the prediction from a transformed input and the transformed prediction from the original input. The same spatial transformation is used in both paths. Comparing these outputs measures equivariance.

\paragraph{Rotation equivariance.} Optical tests rotate the mask together with the quadrature points of the illumination source. Resist tests rotate the input aerial image. We use $90^\circ$, $180^\circ$, and $270^\circ$ rotations, which exactly permute the square grid and require no interpolation. For mask-to-AI prediction, mean equivariance RMSE over the 200 test samples on ViaSet is $3.333\times10^{-9}$, $3.475\times10^{-9}$, and $3.319\times10^{-9}$ at $90^\circ$, $180^\circ$, and $270^\circ$, respectively. MetalSet gives $1.353\times10^{-8}$, $1.462\times10^{-8}$, and $1.343\times10^{-8}$. For AI-to-RI prediction, the respective values are $1.633\times10^{-8}$, $2.106\times10^{-8}$, and $1.641\times10^{-8}$ on ViaSet, and $3.848\times10^{-8}$, $5.354\times10^{-8}$, and $4.040\times10^{-8}$ on MetalSet. The largest absolute residual across the tested cases is $3.159\times10^{-6}$, and the two printed patterns have IOU $1.0000$. These residuals verify rotation equivariance for both tasks at the tested angles, up to floating point error.

\begin{figure}[htbp]
    \centering
    \includegraphics[width=0.99\linewidth]{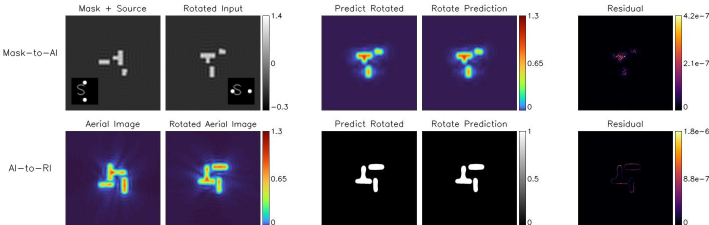}
    \caption{Rotation equivariance of \methodname{} at $90^\circ$ for the two prediction tasks on MetalSet.}
    \label{fig:rotation-equivariance}
\end{figure}

Figure~\ref{fig:rotation-equivariance} compares the output fields for a $90^\circ$ rotation. ``Predict Rotated'' comes from predicting the rotated input. To obtain ``Rotate Prediction'', the model predicts the original input and that output is then rotated by $90^\circ$. The two fields look indistinguishable in each row. Their small differences lie near feature boundaries, consistent with the quantitative residuals.

\paragraph{Translation equivariance.} For mask-to-AI prediction, only the mask is translated; the illumination source stays fixed in pupil coordinates. For AI-to-RI prediction, we periodically translate the input aerial image. The tested shifts are $(0.1,0.1)$, $(0.2,0.2)$, and $(0.5,0.5)$ pixel. They correspond to physical displacements $\boldsymbol{\tau}=(1.6,1.6)\,\mathrm{nm}$, $(3.2,3.2)\,\mathrm{nm}$, and $(8.0,8.0)\,\mathrm{nm}$ for mask-to-AI, and $\boldsymbol{\tau}=(0.6006,0.6006)\,\mathrm{nm}$, $(1.2012,1.2012)\,\mathrm{nm}$, and $(3.0029,3.0029)\,\mathrm{nm}$ for AI-to-RI. To limit nonlinear aliasing on the coarse task grids, we periodically interpolate the input onto a twofold internal grid, halve the runtime pixel spacing, and spectrally project the prediction onto the original grid. Mask-to-AI RMSE values are $9.094\times10^{-9}$, $1.306\times10^{-8}$, and $4.054\times10^{-9}$ on ViaSet, and $3.434\times10^{-6}$, $5.554\times10^{-6}$, and $2.217\times10^{-8}$ on MetalSet. AI-to-RI values are $2.542\times10^{-7}$, $4.108\times10^{-7}$, and $1.520\times10^{-8}$ on ViaSet, and $1.660\times10^{-5}$, $2.677\times10^{-5}$, and $4.251\times10^{-8}$ on MetalSet. The first two shifts remain fractional on the internal grid. The final shift is integer there and serves as a control. \methodname{} thus retains low residuals for subpixel translations, but its discrete implementation does not preserve fractional shifts on the grid to floating point precision.

\begin{figure}[htbp]
    \centering
    \includegraphics[width=0.99\linewidth]{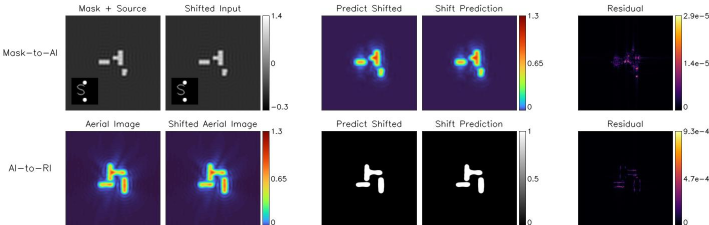}
    \caption{Periodic translation equivariance of \methodname{} for a subpixel shift of $0.1$ pixel on MetalSet.}
    \label{fig:translation-equivariance}
\end{figure}

The $0.1$-pixel translation in Figure~\ref{fig:translation-equivariance} illustrates the subpixel test. ``Predict Shifted'' is predicted from the shifted input. ``Shift Prediction'' is formed by predicting the original input and then shifting the output by $0.1$ pixel in each coordinate direction. Differences occur mainly near feature boundaries. The close field agreement is consistent with the low residuals for subpixel periodic translation.
\begin{table}[H]
    \centering
    \caption{Relative equivariance residuals at a $90^\circ$ rotation and a $0.37$-pixel periodic translation.}
    \label{tab:comparison-equivariance}
    \scriptsize
    \setlength{\tabcolsep}{3pt}
    \begin{tabular}{lllcccc}
        \toprule
        Task & Dataset & Model & Rot. Rel. $L_2\downarrow$ & Rot. IOU $\uparrow$ & Shift Rel. $L_2\downarrow$ & Shift IOU $\uparrow$ \\
        \midrule
        \multirow{8}{*}{Mask$\to$AI}
        & \multirow{4}{*}{ViaSet} & FNO & $3.064{\times}10^{-2}$ & -- & $\mathbf{1.572{\times}10^{-6}}$ & -- \\
        & & DOINN & $5.570{\times}10^{-2}$ & -- & $4.666{\times}10^{-2}$ & -- \\
        & & CFNO & $6.448{\times}10^{-2}$ & -- & $3.631{\times}10^{-2}$ & -- \\
        & & \methodname{} & $\mathbf{2.353{\times}10^{-7}}$ & -- & $\mathbf{7.453{\times}10^{-7}}$ & -- \\
        \cmidrule(lr){2-7}
        & \multirow{4}{*}{MetalSet} & FNO & $2.428{\times}10^{-2}$ & -- & $\mathbf{6.431{\times}10^{-5}}$ & -- \\
        & & DOINN & $5.678{\times}10^{-2}$ & -- & $2.179{\times}10^{-2}$ & -- \\
        & & CFNO & $6.874{\times}10^{-2}$ & -- & $2.247{\times}10^{-2}$ & -- \\
        & & \methodname{} & $\mathbf{1.687{\times}10^{-7}}$ & -- & $\mathbf{5.299{\times}10^{-5}}$ & -- \\
        \midrule
        \multirow{8}{*}{AI$\to$RI}
        & \multirow{4}{*}{ViaSet} & FNO & $3.588{\times}10^{-3}$ & $0.9981$ & $\mathbf{6.268{\times}10^{-6}}$ & $\mathbf{1.0000}$ \\
        & & DOINN & $2.171{\times}10^{-2}$ & $0.9872$ & $1.001{\times}10^{-2}$ & $0.9942$ \\
        & & CFNO & $2.729{\times}10^{-2}$ & $0.9835$ & $7.597{\times}10^{-2}$ & $0.6753$ \\
        & & \methodname{} & $\mathbf{2.797{\times}10^{-7}}$ & $\mathbf{1.0000}$ & $\mathbf{5.396{\times}10^{-6}}$ & $\mathbf{1.0000}$ \\
        \cmidrule(lr){2-7}
        & \multirow{4}{*}{MetalSet} & FNO & $3.906{\times}10^{-3}$ & $0.9984$ & $\mathbf{1.945{\times}10^{-5}}$ & $\mathbf{1.0000}$ \\
        & & DOINN & $2.173{\times}10^{-2}$ & $0.9907$ & $5.655{\times}10^{-3}$ & $0.9976$ \\
        & & CFNO & $2.194{\times}10^{-2}$ & $0.9904$ & $5.161{\times}10^{-3}$ & $0.9978$ \\
        & & \methodname{} & $\mathbf{2.333{\times}10^{-7}}$ & $\mathbf{1.0000}$ & $1.244{\times}10^{-4}$ & $\mathbf{1.0000}$ \\
        \bottomrule
    \end{tabular}
\end{table}

Table~\ref{tab:comparison-equivariance} compares a $90^\circ$ rotation and a subpixel periodic translation of $0.37$ pixel. All four architectures use the same subpixel evaluation protocol. Rotation residuals for \methodname{} are several orders of magnitude below those of the comparison models. Translation residuals are similar to FNO in most cases, while DOINN and CFNO give larger residuals. These observations agree with the construction: radial sharing imposes rotation equivariance, and the Fourier multiplier preserves translation equivariance under periodic boundaries.

\subsection{Transfer from a Single Training Window}
\label{sec:single-window-transfer}

We evaluate direct transfer of the reference-window models under changes in physical window extent and spatial resolution. The reference field is $2048\,\mathrm{nm}\times2048\,\mathrm{nm}$, sampled at the pixel size for each task. Window transfer doubles both side lengths without changing pixel size. Resolution transfer halves pixel size on the same physical window. Joint transfer doubles the side lengths and halves pixel size at once. None of the settings involves retraining, and each target is simulated independently on its evaluation grid.

\begin{table}[H]
    \centering
    \caption{Relative $L_2$ errors under changes in physical window extent and spatial resolution.}
    \label{tab:comparison-size-transfer}
    \scriptsize
    \setlength{\tabcolsep}{5pt}
    \begin{tabular}{lllccc}
        \toprule
        Task & Dataset & Model & Window $\downarrow$ & Resolution $\downarrow$ & Joint $\downarrow$ \\
        \midrule
        \multirow{8}{*}{Mask$\to$AI}
        & \multirow{4}{*}{ViaSet} & FNO & $0.3183$ & $0.0561$ & $0.3271$ \\
        & & DOINN & $0.1147$ & $2.4334$ & $2.4453$ \\
        & & CFNO & $0.1021$ & $1.4297$ & $1.4342$ \\
        & & \methodname{} & $\mathbf{0.0339}$ & $\mathbf{0.0357}$ & $\mathbf{0.0361}$ \\
        \cmidrule(lr){2-6}
        & \multirow{4}{*}{MetalSet} & FNO & $0.6118$ & $0.0311$ & $0.6139$ \\
        & & DOINN & $0.0672$ & $0.6497$ & $0.6481$ \\
        & & CFNO & $0.0973$ & $0.6235$ & $0.6224$ \\
        & & \methodname{} & $\mathbf{0.0250}$ & $\mathbf{0.0275}$ & $\mathbf{0.0276}$ \\
        \midrule
        \multirow{8}{*}{AI$\to$RI}
        & \multirow{4}{*}{ViaSet} & FNO & $0.3793$ & $\mathbf{0.0065}$ & $0.3771$ \\
        & & DOINN & $0.0249$ & $0.1691$ & $0.1691$ \\
        & & CFNO & $0.0310$ & $0.9729$ & $0.9738$ \\
        & & \methodname{} & $\mathbf{0.0182}$ & $0.0131$ & $\mathbf{0.0182}$ \\
        \cmidrule(lr){2-6}
        & \multirow{4}{*}{MetalSet} & FNO & $0.2701$ & $\mathbf{0.0055}$ & $0.2698$ \\
        & & DOINN & $0.0251$ & $0.0622$ & $0.0619$ \\
        & & CFNO & $0.0293$ & $0.1084$ & $0.1063$ \\
        & & \methodname{} & $\mathbf{0.0149}$ & $0.0088$ & $\mathbf{0.0149}$ \\
        \bottomrule
    \end{tabular}
\end{table}

Table~\ref{tab:comparison-size-transfer} reports transfer errors when window extent, resolution, or both change. \methodname{} has the lowest errors in all four window-transfer cases and all four joint-transfer cases. In resolution transfer it gives the lowest errors for mask-to-AI on both datasets, whereas FNO gives the lowest errors for AI-to-RI. Within each task and dataset, its errors stay of the same order across the three settings. This indicates limited dependence on the training window or its discrete resolution. The behavior agrees with evaluating a continuous spectral function on each input's physical frequency grid instead of using a fixed discrete grid. Alongside the prediction and equivariance results, the comparisons show that \methodname{} combines compact parameterization, geometric equivariance, and stable transfer over the evaluated physical windows and resolutions.

\subsection{Prediction on Multiple Windows and Local Error Analysis}
\label{sec:spatial-adaptivity-results}

We next examine spatial adaptivity using dedicated models trained jointly on $2048\,\mathrm{nm}\times2048\,\mathrm{nm}$ and $4096\,\mathrm{nm}\times4096\,\mathrm{nm}$ fields. Each window provides 800 training layouts and 100 validation layouts. Composites use only base layouts from the corresponding split. Within each epoch, batches from the two windows are randomly interleaved. We select the checkpoint with the lowest mean validation relative $L_2$ error across the two training windows. Evaluation covers both seen extents and the unseen $6144\,\mathrm{nm}\times6144\,\mathrm{nm}$ and $6144\,\mathrm{nm}\times4096\,\mathrm{nm}$ fields. For each task, the pixel size stays fixed across all windows. The three larger windows consist of $2\times2$, $3\times3$, and $3\times2$ arrangements of distinct binary layouts, respectively. We first assemble each complete binary field on the original grid, apply the same physical low-pass cutoff, and then resample it onto the task grid by area averaging. Aerial and resist targets are simulated over the complete composed field; they are not assembled from separately simulated subfields. Appendix~\ref{app:training-configuration} summarizes the losses and optimization settings for these dedicated models.

\paragraph{Mask-to-AI prediction.}
\begin{figure}[htbp]
    \centering
    \includegraphics[width=0.96\linewidth]{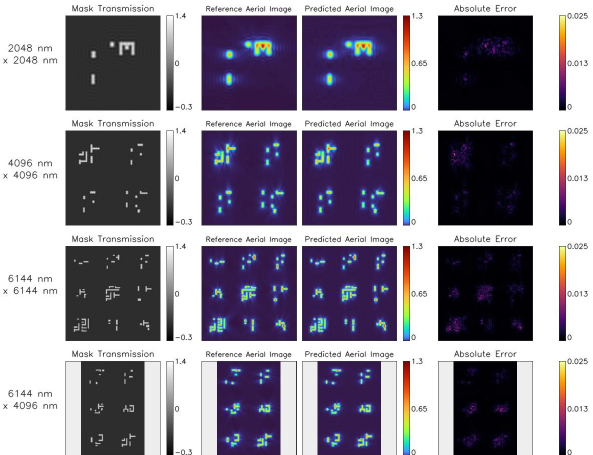}
    \caption{Aerial image predictions for composed MetalSet layouts with different physical window extents.}
    \label{fig:metalset-mask-to-ai-multisize-predictions}
\end{figure}

Figure~\ref{fig:metalset-mask-to-ai-multisize-predictions} presents representative MetalSet predictions for the two seen and two unseen physical windows. At each window, the displayed layout has a maximum absolute error near the median, so these examples describe typical performance rather than the error tail. The predicted aerial images reproduce the overall optical response around the mask features at all four windows. Their maximum absolute errors are $0.0117$, $0.0187$, $0.0223$, and $0.0219$, respectively. Differences are largest near feature boundaries, where aerial intensity changes most rapidly. We examine the distribution of these largest local errors over the full test sets below.

\begin{table}[H]
    \centering
    \small
    \setlength{\tabcolsep}{5pt}
    \caption{Mask-to-AI prediction performance across physical window extents.}
    \label{tab:mask-to-ai-spatial-adaptivity}
    \begin{tabular}{llcc}
        \toprule
        Dataset & Physical window & Output RMSE $\downarrow$ & Output Rel. $L_2$ $\downarrow$ \\
        \midrule
        \multirow{4}{*}{ViaSet}
        & $2048\,\mathrm{nm}\times2048\,\mathrm{nm}$ & $1.853{\times}10^{-4}$ & $0.0131$ \\
        & $4096\,\mathrm{nm}\times4096\,\mathrm{nm}$ & $1.836{\times}10^{-4}$ & $0.0132$ \\
        & $6144\,\mathrm{nm}\times6144\,\mathrm{nm}$ & $1.884{\times}10^{-4}$ & $0.0135$ \\
        & $6144\,\mathrm{nm}\times4096\,\mathrm{nm}$ & $1.863{\times}10^{-4}$ & $0.0132$ \\
        \midrule
        \multirow{4}{*}{MetalSet}
        & $2048\,\mathrm{nm}\times2048\,\mathrm{nm}$ & $1.097{\times}10^{-3}$ & $0.0136$ \\
        & $4096\,\mathrm{nm}\times4096\,\mathrm{nm}$ & $1.086{\times}10^{-3}$ & $0.0136$ \\
        & $6144\,\mathrm{nm}\times6144\,\mathrm{nm}$ & $1.120{\times}10^{-3}$ & $0.0138$ \\
        & $6144\,\mathrm{nm}\times4096\,\mathrm{nm}$ & $1.142{\times}10^{-3}$ & $0.0140$ \\
        \bottomrule
    \end{tabular}
\end{table}

We then test all 200 layouts from each dataset at every physical window to assess whether the representative behavior holds more broadly. Table~\ref{tab:mask-to-ai-spatial-adaptivity} lists aerial-image RMSE and relative $L_2$ error. Across the four windows, RMSE lies between $1.836\times10^{-4}$ and $1.884\times10^{-4}$ on ViaSet, and between $1.086\times10^{-3}$ and $1.142\times10^{-3}$ on MetalSet. Relative $L_2$ errors stay within $0.0131$--$0.0135$ and $0.0136$--$0.0140$, respectively. Neither measure increases systematically with physical window area. Thus, the same model retains nearly constant field averaged aerial-image accuracy on the larger windows, including both unseen extents.

\paragraph{AI-to-RI prediction.}
\begin{figure}[htbp]
    \centering
    \includegraphics[width=0.96\linewidth]{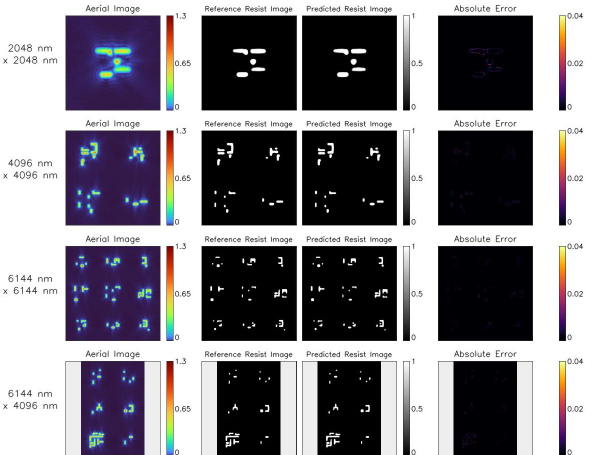}
    \caption{Resist image predictions for composed MetalSet layouts with different physical window extents.}
    \label{fig:metalset-multisize-predictions}
\end{figure}

The corresponding MetalSet AI-to-RI predictions appear in Figure~\ref{fig:metalset-multisize-predictions}. The four displayed windows have maximum absolute errors of $0.0210$, $0.0282$, $0.0350$, and $0.0321$, respectively. Differences occur mainly within narrow bands around resist boundaries. Although locally noticeable, these differences have little effect after thresholding. The predicted printed patterns consequently remain visually consistent with the simulated patterns at all four physical windows.

\begin{table}[H]
    \centering
    \small
    \setlength{\tabcolsep}{4pt}
    \caption{AI-to-RI prediction performance across physical window extents.}
    \label{tab:ai-to-ri-spatial-adaptivity}
    \begin{tabular}{llcccc}
        \toprule
        Dataset & Physical window & Output RMSE $\downarrow$ & Output Rel. $L_2$ $\downarrow$ & IOU $\uparrow$ & PA $\uparrow$ \\
        \midrule
        \multirow{4}{*}{ViaSet}
        & $2048\,\mathrm{nm}\times2048\,\mathrm{nm}$ & $4.780{\times}10^{-4}$ & $0.0082$ & $0.9954$ & $0.9969$ \\
        & $4096\,\mathrm{nm}\times4096\,\mathrm{nm}$ & $5.926{\times}10^{-4}$ & $0.0102$ & $0.9941$ & $0.9969$ \\
        & $6144\,\mathrm{nm}\times6144\,\mathrm{nm}$ & $5.740{\times}10^{-4}$ & $0.0099$ & $0.9943$ & $0.9970$ \\
        & $6144\,\mathrm{nm}\times4096\,\mathrm{nm}$ & $5.874{\times}10^{-4}$ & $0.0101$ & $0.9941$ & $0.9969$ \\
        \midrule
        \multirow{4}{*}{MetalSet}
        & $2048\,\mathrm{nm}\times2048\,\mathrm{nm}$ & $6.701{\times}10^{-4}$ & $0.0041$ & $0.9982$ & $0.9984$ \\
        & $4096\,\mathrm{nm}\times4096\,\mathrm{nm}$ & $7.185{\times}10^{-4}$ & $0.0044$ & $0.9980$ & $0.9984$ \\
        & $6144\,\mathrm{nm}\times6144\,\mathrm{nm}$ & $7.406{\times}10^{-4}$ & $0.0045$ & $0.9980$ & $0.9985$ \\
        & $6144\,\mathrm{nm}\times4096\,\mathrm{nm}$ & $7.539{\times}10^{-4}$ & $0.0045$ & $0.9979$ & $0.9984$ \\
        \bottomrule
    \end{tabular}
\end{table}

Table~\ref{tab:ai-to-ri-spatial-adaptivity} reports the AI-to-RI results for all 200 test layouts at each window. On ViaSet, RMSE rises from $4.780\times10^{-4}$ at the reference window to $5.740\times10^{-4}$--$5.926\times10^{-4}$ at the three larger windows. Relative $L_2$ error similarly rises from $0.0082$ to $0.0099$--$0.0102$, with comparable errors across the larger windows. The changes on MetalSet are more gradual: RMSE increases from $6.701\times10^{-4}$ to $7.539\times10^{-4}$, and relative $L_2$ error changes from $0.0041$ to $0.0045$. Even with these changes in the continuous resist images, the reported IOU values are at least $0.9940$ on ViaSet and $0.9979$ on MetalSet, while PA is at least $0.9968$ and $0.9984$. Printed patterns therefore remain stable at both seen and unseen extents. The moderate growth in continuous field error is not accompanied by a comparable loss of agreement in the printed patterns.

\paragraph{Distribution of localized errors.}
RMSE and relative $L_2$ error measure discrepancies averaged over the whole field. By contrast, the maximum absolute error of each layout captures isolated local differences. For every task and dataset, we set a reference error level on the $2048\,\mathrm{nm}\times2048\,\mathrm{nm}$ window such that $98.0\%$ of its 200 test layouts, or 196 layouts, fall at or below that level. We then use this fixed level for the remaining windows. For mask-to-AI prediction, the levels are $0.009$ on ViaSet and $0.031$ on MetalSet; for AI-to-RI prediction, they are $0.057$ and $0.043$, respectively. Coverage is the fraction of layouts whose maximum error falls at or below the corresponding fixed level. The reference-window value of $98.0\%$ is therefore a calibration baseline; the comparison concerns how coverage changes across physical windows. We calibrate each task and dataset separately because their error distributions differ. These levels serve as references for window transfer, rather than common manufacturing acceptance limits.

\begin{figure}[htbp]
    \centering
    \includegraphics[width=\linewidth,clip]{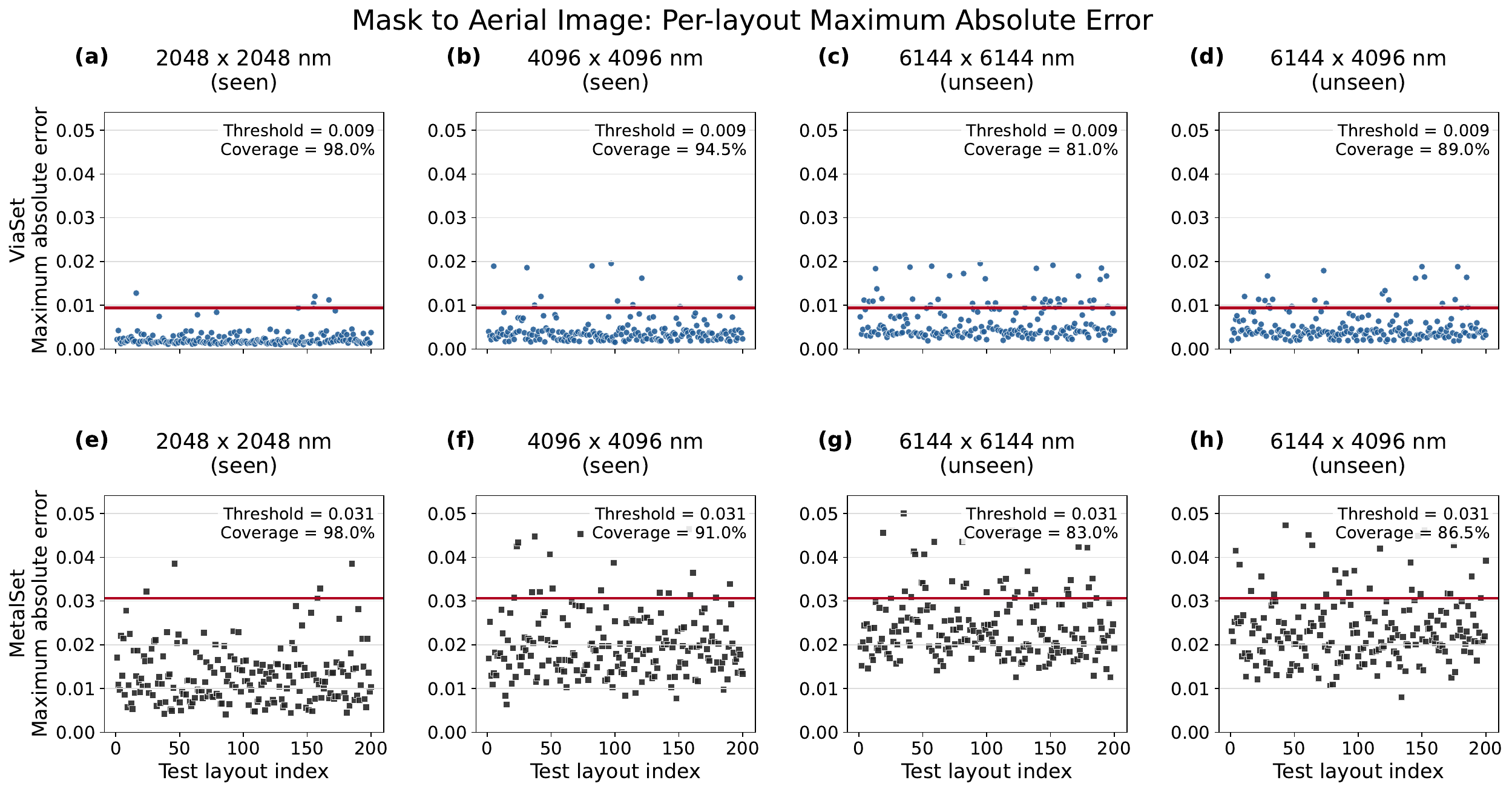}
    \caption{Per-layout maximum absolute errors for mask-to-AI prediction across physical window extents. The red lines denote the fixed reference error levels defined in the text.}
    \label{fig:mask-to-ai-multisize-error-distribution}
\end{figure}

Figure~\ref{fig:mask-to-ai-multisize-error-distribution} reports the mask-to-AI distributions. Starting from $98.0\%$, ViaSet coverage falls to $94.5\%$, $81.0\%$, and $89.0\%$ on the $4096\,\mathrm{nm}\times4096\,\mathrm{nm}$, $6144\,\mathrm{nm}\times6144\,\mathrm{nm}$, and $6144\,\mathrm{nm}\times4096\,\mathrm{nm}$ windows, respectively. MetalSet coverage across the four windows is $98.0\%$, $91.0\%$, $83.0\%$, and $86.5\%$. In both datasets, coverage is lowest at the largest $6144\,\mathrm{nm}\times6144\,\mathrm{nm}$ window and recovers partly at the smaller $6144\,\mathrm{nm}\times4096\,\mathrm{nm}$ window. With a fixed pixel size, these windows contain approximately $1$, $4$, $9$, and $6$ times as many sampled locations. A maximum taken over a larger field thus has more opportunities to encounter an isolated large error. This also accounts for the lower coverage at the seen $4096\,\mathrm{nm}\times4096\,\mathrm{nm}$ window than at the reference window: coverage describes the tail of localized errors, not simply whether an extent was seen during training. Meanwhile, Table~\ref{tab:mask-to-ai-spatial-adaptivity} shows almost unchanged RMSE and relative $L_2$ errors. Together, these observations indicate broader localized error tails rather than a general loss of aerial-image accuracy.

\begin{figure}[htbp]
    \centering
    \includegraphics[width=\linewidth,clip]{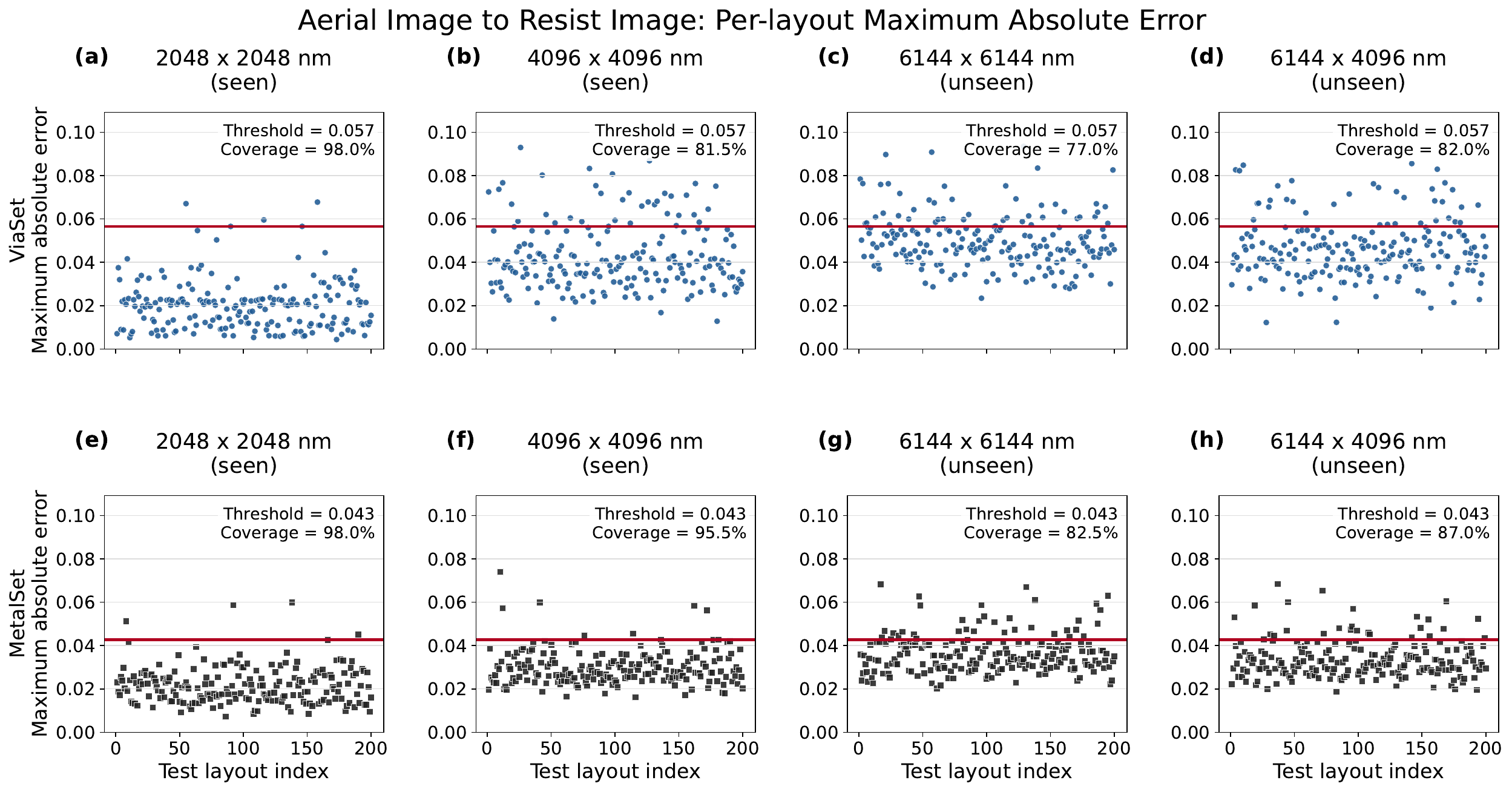}
    \captionsetup{width=0.96\linewidth}
    \caption{Per-layout maximum absolute errors for AI-to-RI prediction across physical window extents. The red lines denote the fixed reference error levels defined in the text.}
    \label{fig:ai-to-ri-multisize-error-distribution}
\end{figure}

The AI-to-RI distributions are shown in Figure~\ref{fig:ai-to-ri-multisize-error-distribution}. Coverage on ViaSet is $98.0\%$, $81.5\%$, $77.0\%$, and $82.0\%$ across the four windows; on MetalSet, it is $98.0\%$, $95.5\%$, $82.5\%$, and $87.0\%$. Both datasets again have their lowest coverage at $6144\,\mathrm{nm}\times6144\,\mathrm{nm}$. Reducing the area to $6144\,\mathrm{nm}\times4096\,\mathrm{nm}$ raises coverage by $5.0$ percentage points on ViaSet and $4.5$ on MetalSet. The ViaSet values of $81.5\%$ and $82.0\%$ at $4096\,\mathrm{nm}\times4096\,\mathrm{nm}$ and $6144\,\mathrm{nm}\times4096\,\mathrm{nm}$ differ by just one of the 200 test layouts. This small difference does not establish a systematic reversal of the area trend. In resist prediction, a slight displacement of a steep transition can cause a large pointwise difference within a narrow band. The maximum error can therefore cross the reference level even when the resist image is accurate over most of the field. This agrees with Table~\ref{tab:ai-to-ri-spatial-adaptivity}: field averaged errors change moderately, while IOU and PA remain close to one at all four windows.

\section{Discussion and Conclusions}
\label{sec:conclusion}

We have presented \methodname{}, a spatially adaptive and geometrically equivariant neural operator for lithography simulation. Its basic spectral parameterization represents Fourier weights as continuous functions of normalized radial frequency through a learnable sum of Gaussian functions. Radial sharing enforces rotation equivariance in the continuous construction, while Fourier multipliers preserve translation equivariance under periodic boundary conditions. Evaluating the continuous spectral response on the physical frequency grids of different layout windows allows the learned parameters to be applied across windows of different physical extents. The localized Gaussian basis and smooth cutoff taper regularize the spectral response and promote spatial decay of the induced nonlocal kernels.

The main experiments demonstrate geometric equivariance and spatial adaptivity in both prediction tasks. Rotation residuals at $90^\circ$, $180^\circ$, and $270^\circ$ remain near floating point precision, with the mask and illumination source rotated jointly in the optical task. Subpixel periodic translations also yield low residuals and closely aligned output fields. Models trained jointly on $2048\,\mathrm{nm}\times2048\,\mathrm{nm}$ and $4096\,\mathrm{nm}\times4096\,\mathrm{nm}$ fields apply directly to unseen window extents and shapes, without retraining or resizing inputs to the training window. Across the evaluated windows, errors in aerial image prediction remain almost unchanged. errors in resist image prediction increase moderately before stabilizing at the larger windows, while printed patterns maintain high IOU and PA. The distributions of the maximum error for individual layouts develop broader upper tails as the evaluated area grows, and the largest square window has the lowest reference coverage. Field averaged metrics remain much more stable. These results demonstrate accurate field predictions and high agreement between the predicted and reference printed patterns across the evaluated physical windows.

Comparisons with FNO, DOINN, and CFNO further assess predictive accuracy and parameter efficiency alongside these structural properties. Rotation equivariance residuals for \methodname{} are substantially lower than those of the comparison models, demonstrating stronger consistency under rotation. For periodic translation, its residuals are comparable to the lowest values in most cases, while the predicted printed patterns maintain nearly perfect overlap with the translated reference patterns on both datasets. Models trained on the reference window also maintain low prediction errors when physical window extent and spatial resolution change, achieving the lowest errors among the compared models in the window transfer and joint transfer tests. On the reference window, \methodname{} achieves competitive predictive accuracy with substantially fewer parameters. These results demonstrate that \methodname{} combines accurate prediction, rotation and translation equivariance, and spatial adaptivity in a compact lithography surrogate. Future work will improve localized pointwise accuracy, explore more flexible equivariant spectral representations, and extend evaluation to a wider range of optical and resist process conditions.

\section*{Acknowledgments}
\addcontentsline{toc}{section}{Acknowledgments}

{\emergencystretch=3em
This work is financially supported by National Key R\&D Program of China (grant No.2024YFA1012403) and the National Natural Science Foundation of China (grant Nos.12325113 and 12426304). The authors thank Dr. Yihao Liang for his helpful suggestions to this work. The authors also acknowledge support from the SJTU Kunpeng \& Ascend Center of Excellence.
\par}

\bibliographystyle{lithosage-compact}
\bibliography{references}

\appendix
\numberwithin{equation}{section}
\addcontentsline{toc}{section}{Appendix}
\section{Configuration and Notation}
\label{app:configuration}

This appendix lists the numerical parameters and implementation settings for the experiments. The main text defines the symbols needed to follow the physical models, data pipeline, and network construction. The following tables give their values and distinguish related notation that might otherwise cause confusion.

\subsection{Notation and Parameter Conventions}
\label{app:notation-configuration}

\begin{table}[H]
    \centering
    \footnotesize
    \setlength{\tabcolsep}{4pt}
    \renewcommand{\arraystretch}{1.10}
    \caption{Key notation and parameter conventions used across the physical model, network, and experiments.}
    \label{tab:notation-conventions}
    \begin{tabularx}{\linewidth}{@{}P{0.15\linewidth}P{0.28\linewidth}Y@{}}
        \toprule
        Symbol & Meaning & Convention or distinction \\
        \midrule
        $\chi,m,a,z,p$ & Layout and lithography fields & Binary layout, mask transmission, aerial image, resist image, and printed pattern, respectively. \\
        $S(\boldsymbol{s})$; $S_Q$ & Source representations & $S$ is continuous; $S_Q=\{(\boldsymbol{s}_q,w_q)\}_{q=1}^{N_s}$ is the quadrature representation used in computation. \\
        $\boldsymbol{r}=(x,y)$ & Lateral spatial coordinate & Measured in $\mathrm{nm}$; $r_{\mathrm{sep}}$ denotes scalar spatial separation. \\
        $\boldsymbol{f}$ & Spatial frequency & Cycles per unit length; used in the optical model and mask preprocessing. \\
        $\boldsymbol{\xi}=(\xi_x,\xi_y)$ & Angular spatial frequency & $\boldsymbol{\xi}=2\pi\boldsymbol{f}$; used in the spectral parameterization of \methodname{}. \\
        $f_c$ & Coherent optical cutoff & Optical pupil cutoff, $f_c=\mathrm{NA}/\lambda$. \\
        $f_{\mathrm{LP}}$ & Mask preprocessing cutoff & Applied before optical simulation; distinct from the network cutoff $\xi_c$. \\
        $\xi_c$ & \methodname{} spectral cutoff & Angular frequency cutoff that defines the retained mode set. \\
        $L_{\mathrm{ref}}$ & Reference side length & Side length of the square reference window. \\
        $L$ & Number of Fourier layers & Network depth; distinct from $L_x$, $L_y$, and $L_{\mathrm{ref}}$. \\
        $d$ & Latent channel width & Upright $\mathrm{d}$ is reserved for integration differentials. \\
        $N$, $N_s$, $N_b$ & Pixel, source point, and basis counts & The subscripts distinguish the three quantities. \\
        \bottomrule
    \end{tabularx}
\end{table}

\subsection{Physical Simulation and Data Generation Parameters}
\label{app:physical-data-configuration}

\begin{table}[H]
    \centering
    \footnotesize
    \setlength{\tabcolsep}{4pt}
    \renewcommand{\arraystretch}{1.10}
    \caption{Spatial discretization and optical data generation parameters.}
    \label{tab:physical-data-parameters}
    \begin{tabularx}{\linewidth}{@{}P{0.28\linewidth}P{0.19\linewidth}Y@{}}
        \toprule
        Parameter & Value & Role \\
        \midrule
        Reference side length $L_{\mathrm{ref}}$ & $2048\,\mathrm{nm}$ & Reference physical window. \\
        Grid of the original binary mask & $2048\times2048$ & $1\,\mathrm{nm}$ sampling before filtering and resampling. \\
        Mask low-pass cutoff $f_{\mathrm{LP}}$ & $0.0127\,\mathrm{nm}^{-1}$ & Approximate transmitted frequency support over all source points. \\
        Mask-to-AI grid & $128\times128$ & Pixel size $16\,\mathrm{nm}$. \\
        AI-to-RI grid & $341\times341$ & Pixel size $6.0059\,\mathrm{nm}$. \\
        Wavelength $\lambda$ & $193\,\mathrm{nm}$ & Optical wavelength. \\
        Numerical aperture $\mathrm{NA}$ & $1.35$ & Projection numerical aperture. \\
        Coherent cutoff $f_c$ & $0.0070\,\mathrm{nm}^{-1}$ & $f_c=\mathrm{NA}/\lambda$. \\
        Projection pupil & Ideal circular pupil & No aberration or defocus. \\
        Source distribution $S$ & Horizontal/vertical dipole & Equal representation of both orientations in each split. \\
        Source quadrature size $N_s$ & $32$ & Quadrature points per source. \\
        \bottomrule
    \end{tabularx}
\end{table}

\begin{table}[H]
    \centering
    \footnotesize
    \setlength{\tabcolsep}{4pt}
    \renewcommand{\arraystretch}{1.10}
    \caption{CAR compact model parameters used to generate the AI-to-RI datasets.}
    \label{tab:car-data-parameters}
    \begin{tabularx}{\linewidth}{@{}P{0.23\linewidth}P{0.24\linewidth}Y@{}}
        \toprule
        Parameter & Value & Meaning \\
        \midrule
        $\alpha_{\mathrm{exp}}$ & $2.5$ & Exposure sensitivity. \\
        $D$ & $0.9$ (ViaSet); $0.5$ (MetalSet) & Normalized exposure dose. \\
        $\sigma$ & $16\,\mathrm{nm}$ & Lateral diffusion length. \\
        $k_{\mathrm{amp}}$ & $3$ & Amplification strength. \\
        $n_{\mathrm{amp}}$ & $1$ & Amplification exponent. \\
        $\vartheta_{\mathrm{dev}}$ & $0.5$ & Development threshold. \\
        $\beta$ & $40$ & Smooth threshold sharpness. \\
        Printed pattern threshold & $0.5$ & Threshold for the continuous resist image. \\
        \bottomrule
    \end{tabularx}
\end{table}

The two datasets share all CAR material parameters but use different normalized exposure doses $D$. At the AI-to-RI pixel spacing, $\sigma=16\,\mathrm{nm}$ corresponds to $2.6641$ pixels on the $341\times341$ grid. For each dataset, we choose $D$ from a fixed candidate sweep using the first 100 training layouts only. Selection maximizes aggregate IOU between simulated printed patterns and the original binary layouts. This dose is then held fixed throughout the training, validation, and test sets.

\subsection{\methodname{} Architecture Parameters}
\label{app:lithosage-configuration}

\begin{table}[H]
    \centering
    \footnotesize
    \setlength{\tabcolsep}{4pt}
    \renewcommand{\arraystretch}{1.10}
    \caption{Architecture parameters of \methodname{} for the two prediction tasks.}
    \label{tab:lithosage-parameters}
    \begin{tabularx}{\linewidth}{@{}P{0.20\linewidth}P{0.20\linewidth}P{0.20\linewidth}Y@{}}
        \toprule
        Parameter & Mask-to-AI & AI-to-RI & Role \\
        \midrule
        Fourier layers $L$ & $4$ & $4$ & Spectral layer count. \\
        Latent width $d$ & $32$ & $32$ & Latent channel count. \\
        Reference side length $L_{\mathrm{ref}}$ & $2048\,\mathrm{nm}$ & $2048\,\mathrm{nm}$ & Reference length for the cutoff. \\
        Spectral cutoff $\xi_c$ & $0.0706\,\mathrm{rad}\,\mathrm{nm}^{-1}$ & $0.0706\,\mathrm{rad}\,\mathrm{nm}^{-1}$ & Fixed cutoff; 24 modes on the reference window, with indices recomputed for each window. \\
        Gaussian bases $N_b$ & $16$ & $16$ & Radial basis count. \\
        Gaussian centers $\mu_b$ & $\sqrt{(b-1/2)/N_b}$ & Linearly spaced on $[0,1]$ & Centers in normalized radial frequency. \\
        Gaussian widths $\delta_{\ell,b}$ & Fixed from adjacent center spacing & Learned independently by layer & Basis widths. \\
        Source scale $\alpha_s$ & $0.6228$ & -- & Scale from source coordinates to normalized angular frequency. \\
        \bottomrule
    \end{tabularx}
\end{table}

For mask-to-AI prediction, source conditioning is implemented on a centered full spectrum so that the mask and source can be rotated jointly, with $\alpha_s=2\pi f_c/\xi_c$. The spectral coefficients use an initialization gain of $0.5$. The final weight matrix in the output head uses a gain of $0.1$. Its output bias is initialized from four times the mean target intensity in the training split. For AI-to-RI prediction, we use a real radial spectral multiplier, in accordance with the real-field symmetry described in Appendix~\ref{app:real-spectral-weights}. The spectral response is multiplied by the cosine taper
\begin{equation}
    t(\nu)
    =
    \begin{cases}
        1,
        & 0\leq\nu\leq0.9,\\[0.3em]
        \dfrac{1}{2}\left[
        1+\cos\!\left(
        \pi\dfrac{\nu-0.9}{0.1}
        \right)
        \right],
        & 0.9<\nu<1,\\[0.8em]
        0,
        & \nu\geq1,
    \end{cases}
    \label{eq:appendix-cosine-taper}
\end{equation}
in which $\nu=\sqrt{\xi_x^2+\xi_y^2}/\xi_c$.

\subsection{Training and Evaluation Parameters}
\label{app:training-configuration}

\begin{table}[H]
    \centering
    \footnotesize
    \setlength{\tabcolsep}{4pt}
    \renewcommand{\arraystretch}{1.10}
    \caption{Training settings for \methodname{} in the reference window and multiple window experiments.}
    \label{tab:training-parameters}
    \begin{tabularx}{\linewidth}{@{}P{0.20\linewidth}P{0.31\linewidth}Y@{}}
        \toprule
        Setting & Reference window & Multiple windows \\
        \midrule
        Training windows & $2048\,\mathrm{nm}\times2048\,\mathrm{nm}$ & $2048\,\mathrm{nm}\times2048\,\mathrm{nm}$ and $4096\,\mathrm{nm}\times4096\,\mathrm{nm}$ \\
        Batch size & $2$ & $2$ for each window loader \\
        Random seed & $0$ & $0$ \\
        Optimizer & Adam, $(\beta_1,\beta_2)=(0.9,0.999)$; no weight decay & Same \\
        Rotation augmentation & $0^\circ,90^\circ,180^\circ,270^\circ$ sampled uniformly & Same \\
        Maximum epochs & $20$ & $20$ \\
        Mask-to-AI loss & Pixelwise MSE & Pixelwise MSE plus $0.2$ times the MSE of the largest $1\%$ absolute residuals \\
        AI-to-RI loss & Pixelwise MSE & Pixelwise MSE \\
        Checkpoint criterion & Lowest validation relative $L_2$ error & Lowest mean validation relative $L_2$ error over both windows \\
        \bottomrule
    \end{tabularx}
\end{table}

For \methodname{}, mask-to-AI training uses a learning rate of $2\times10^{-3}$ in epochs 1--10 and $2\times10^{-4}$ in epochs 11--20. AI-to-RI training uses $10^{-3}$ in epochs 1--5 and $10^{-4}$ in epochs 6--20. Original aerial images, mask-to-AI targets, and AI-to-RI inputs are scaled to $[0,2]$; continuous RI targets and AI-to-RI outputs are in $[0,1]$. The optical output is twice a sigmoid output, whereas the RI output uses a sigmoid directly. Neither training nor checkpoint selection uses the test split.

The experiment subsections give the equivariance settings and localized-error evaluation parameters, since these specify the reported tests directly. In particular, Sec.~\ref{sec:equivariance-results} defines the subpixel translations, and Sec.~\ref{sec:spatial-adaptivity-results} defines the reference levels for the localized maximum-error analysis.

\subsection{Baseline Configurations}
\label{app:baseline-configuration}

Although the three comparison architectures have different internal representations, none natively accepts an unordered set of source quadrature points. We therefore supply the illumination condition to FNO, DOINN, and CFNO through a common source descriptor for mask-to-AI prediction. For a function $\psi$ of source coordinates,
\begin{equation}
    \mathbb{E}_{w}[\psi(\boldsymbol{s})]
    =
    \sum_{q=1}^{N_s}w_q\psi(\boldsymbol{s}_q),
    \label{eq:appendix-weighted-expectation}
\end{equation}
with $\sum_qw_q=1$. We use the following nine moments in the descriptor:
\begin{equation}
    \bigl(
    \mathbb{E}_{w}[s_x],\mathbb{E}_{w}[s_y],
    \mathbb{E}_{w}[s_x^2],\mathbb{E}_{w}[s_y^2],
    \mathbb{E}_{w}[s_xs_y],\mathbb{E}_{w}[s_x^2+s_y^2],
    \mathbb{E}_{w}[s_x^4],\mathbb{E}_{w}[s_y^4],
    \mathbb{E}_{w}[s_x^2s_y^2]
    \bigr),
    \label{eq:appendix-source-moments}
\end{equation}
A $9\rightarrow32\rightarrow16$ multilayer perceptron maps these moments to a conditioning vector. The moments distinguish the two fixed orthogonal dipoles used here without ambiguity. This encoding gives the baselines a compatible source-conditioning interface. By comparison, \methodname{} incorporates the quadrature set $S_Q$ directly into its spectral layer.

\paragraph{FNO.} The FNO baseline uses the standard spectral architecture described in Sec.~\ref{sec:fourier-framework}. After lifting, a stack of Fourier layers processes the latent field using independent complex weights on the retained modes. For comparison, we match its layer count and latent width to \methodname{}. Its spectral weights are still attached to discrete Fourier modes, rather than continuous radial functions. In mask-to-AI prediction, the source descriptor conditions each spectral block.

\paragraph{DOINN.} DOINN joins a global RFNO branch with a local convolutional path. Pooling supplies the RFNO input, while a local encoder retains higher-frequency spatial information at several channel scales. In mask-to-AI prediction, the source descriptor conditions the RFNO feature and the features at every encoder scale before the two branches are combined.

\paragraph{CFNO.} CFNO extracts spectral features at multiple patch scales and combines them with a convolutional encoder and decoder. Its three patch branches represent different spatial scales, while the convolutional path gathers local features. For mask-to-AI prediction, source conditioning is applied to the spectral and convolutional feature groups before they are combined.

The numerical architecture settings for these implementations are listed in Table~\ref{tab:baseline-parameters}.
\begin{table}[H]
    \centering
    \footnotesize
    \setlength{\tabcolsep}{4pt}
    \renewcommand{\arraystretch}{1.10}
    \caption{Core numerical architecture settings of the comparison models.}
    \label{tab:baseline-parameters}
    \begin{tabularx}{\linewidth}{@{}P{0.12\linewidth}P{0.37\linewidth}Y@{}}
        \toprule
        Model & Spectral or multiscale settings & Local or channel settings \\
        \midrule
        FNO & $4$ Fourier layers; $24$ retained modes & Latent width $32$ \\
        DOINN & $8\times8$ average pooling; RFNO with $24$ retained modes & RFNO output $96$; encoder widths $(24,48,96)$ \\
        CFNO & Patch sizes $(16,32,64)$; feature dimensions $(16,32,64)$ & Encoder $(32,64,128)$; decoder $(128,64,32)$ \\
        \bottomrule
    \end{tabularx}
\end{table}

Table~\ref{tab:baseline-learning-rates} lists each model's learning rate schedule. The remaining shared training settings follow Sec.~\ref{sec:experimental-setup} and Appendix~\ref{app:training-configuration}.
\begin{table}[H]
    \centering
    \footnotesize
    \setlength{\tabcolsep}{4pt}
    \renewcommand{\arraystretch}{1.10}
    \caption{Learning rate schedules for the baseline models trained on the reference window.}
    \label{tab:baseline-learning-rates}
    \begin{tabularx}{\linewidth}{@{}P{0.12\linewidth}Y Y@{}}
        \toprule
        Model & Mask-to-AI & AI-to-RI \\
        \midrule
        FNO & $2\times10^{-3}$ (1--3); $2\times10^{-4}$ (4--10); $10^{-4}$ (11--15); $5\times10^{-5}$ (16--20) & $5\times10^{-4}$ (1--5); $5\times10^{-5}$ (6--10); $2.5\times10^{-5}$ (11--20) \\
        DOINN & $5\times10^{-4}$ (1--5); $5\times10^{-5}$ (6--20) & $5\times10^{-4}$ (1--5); $5\times10^{-5}$ (6--10); $2.5\times10^{-5}$ (11--20) \\
        CFNO & $5\times10^{-4}$ (1--5); $5\times10^{-5}$ (6--10); $2.5\times10^{-5}$ (11--20) & $3\times10^{-4}$ (1--3); $3\times10^{-5}$ (4--10); $2.5\times10^{-5}$ (11--20) \\
        \bottomrule
    \end{tabularx}
\end{table}

\section{Proof of the Radial Spectral Weight Characterization}
\label{app:radial-proof}

\begin{proof}[Proof of Theorem~\ref{thm:radial-equivariance}]
Write $\mathcal{L}_G$ for the Fourier layer in Eq.~\eqref{eq:continuous-fourier-layer}. A rotation through $\alpha\in[0,2\pi)$ acts on a field as
\begin{equation}
    (R_{\alpha}u)(\boldsymbol{x})
    =
    u(R_{\alpha}^{-1}\boldsymbol{x}).
    \label{eq:rotation-action-proof}
\end{equation}
For multiple scalar channels, the rotation acts on each channel separately without mixing channel indices. Changing variables yields
\begin{equation}
    \widehat{R_{\alpha}u}(\boldsymbol{\xi})
    =
    \widehat{u}(R_{\alpha}^{-1}\boldsymbol{\xi}),
    \label{eq:rotated-fourier-transform}
\end{equation}
with $\boldsymbol{\xi}=(\xi_x,\xi_y)$. If we rotate the input before applying the Fourier layer, we obtain
\begin{equation}
    \widehat{\mathcal{L}_G R_{\alpha}u}(\boldsymbol{\xi})
    =
    G(\boldsymbol{\xi})
    \widehat{u}(R_{\alpha}^{-1}\boldsymbol{\xi}).
    \label{eq:operator-after-rotation}
\end{equation}
If we apply the layer before rotating the output, we obtain
\begin{equation}
    \widehat{R_{\alpha}\mathcal{L}_G u}(\boldsymbol{\xi})
    =
    G(R_{\alpha}^{-1}\boldsymbol{\xi})
    \widehat{u}(R_{\alpha}^{-1}\boldsymbol{\xi}).
    \label{eq:rotation-after-operator}
\end{equation}
The layer is therefore rotation equivariant provided that
\begin{equation}
    G(\boldsymbol{\xi})
    =
    G(R_{\alpha}^{-1}\boldsymbol{\xi})
    \label{eq:rotation-invariant-multiplier}
\end{equation}
at every $\boldsymbol{\xi}$ and for every $\alpha$. Direct substitution establishes sufficiency.

For necessity, assume that the layer is rotation equivariant, yet the matrices in Eq.~\eqref{eq:rotation-invariant-multiplier} differ at a frequency $\boldsymbol{\xi}_0$. By continuity of their difference, there exist a channel vector $c$ and a neighborhood $U$ of $\boldsymbol{\xi}_0$ such that
\begin{equation}
    \left[
    G(\boldsymbol{\xi})
    -G(R_{\alpha}^{-1}\boldsymbol{\xi})
    \right]c
    \neq0.
    \label{eq:necessity-matrix-difference}
\end{equation}
Now choose Fourier data supported in $R_{\alpha}^{-1}U$ and aligned with $c$. The outputs in Eqs.~\eqref{eq:operator-after-rotation} and \eqref{eq:rotation-after-operator} then differ on $U$, contradicting rotation equivariance. Thus, Eq.~\eqref{eq:rotation-invariant-multiplier} must hold at every frequency.

A rotation leaves the frequency radius unchanged:
\begin{equation}
    \left\|R_{\alpha}^{-1}\boldsymbol{\xi}\right\|_2
    =
    \left\|\boldsymbol{\xi}\right\|_2.
    \label{eq:rotation-preserves-frequency-radius}
\end{equation}
When $G$ depends only on $\|\boldsymbol{\xi}\|_2$, this identity ensures Eq.~\eqref{eq:rotation-invariant-multiplier} and hence rotation equivariance. Conversely, Eq.~\eqref{eq:rotation-invariant-multiplier} requires equal values of $G$ at frequencies connected by a rotation. Any two points on a circle centered at the origin can be related by a planar rotation. Therefore, $G$ is constant along each such circle and is a function only of $\|\boldsymbol{\xi}\|_2$.
\end{proof}

\section{Real Spectral Weights for Real Fields}
\label{app:real-spectral-weights}

We use the angular spatial frequencies $\xi_x$ and $\xi_y$ defined in Eq.~\eqref{eq:physical-frequency}. Under this convention, the continuous Fourier transform takes the form
\begin{equation}
    \widehat{u}(\xi_x,\xi_y)
    =
    \int_{\mathbb{R}^2}
    u(x,y)
    \exp\!\left[
        -\mathrm{i}(\xi_xx+\xi_yy)
    \right]
    \,\mathrm{d}x\,\mathrm{d}y.
    \label{eq:fourier-transform-real-field}
\end{equation}
If $u$ is a real field, then
\begin{equation}
    \widehat{u}(-\xi_x,-\xi_y)
    =
    \overline{\widehat{u}(\xi_x,\xi_y)}.
    \label{eq:input-conjugate-symmetry}
\end{equation}
The output $y$ must obey the same relation to be real. Together with Eq.~\eqref{eq:continuous-fourier-layer}, this gives the necessary and sufficient condition
\begin{equation}
    G(-\xi_x,-\xi_y)
    =
    \overline{G(\xi_x,\xi_y)}.
    \label{eq:weight-conjugate-symmetry}
\end{equation}
For a matrix-valued multiplier, this equality applies to each entry.

The AI-to-RI spectral weight is radial about the origin, which implies
\begin{equation}
    G(-\xi_x,-\xi_y)
    =
    G(\xi_x,\xi_y).
    \label{eq:opposite-frequency-radial-weight}
\end{equation}
Combining Eqs.~\eqref{eq:weight-conjugate-symmetry} and \eqref{eq:opposite-frequency-radial-weight} yields
\begin{equation}
    G(\xi_x,\xi_y)
    =
    \overline{G(\xi_x,\xi_y)}.
    \label{eq:real-radial-spectral-weight}
\end{equation}
The radial multiplier must consequently be real. The Gaussian basis functions are real, so we enforce this requirement through
\begin{equation}
    \gamma_{\ell,o,i,b}\in\mathbb{R}.
    \label{eq:real-gaussian-coefficients}
\end{equation}

In mask-to-AI prediction, the Gaussian responses are centered relative to source points, so the multiplier is not radial about the origin. The dipole quadrature used here is, however, centrally symmetric. Each source point $(\boldsymbol{s}_q,w_q)$ has a partner $(\boldsymbol{s}_{q^{\ast}},w_{q^{\ast}})$ for which
\begin{equation}
    \boldsymbol{s}_{q^{\ast}}=-\boldsymbol{s}_q,
    \qquad
    w_{q^{\ast}}=w_q.
    \label{eq:source-pair-symmetry}
\end{equation}
Each such pair satisfies
\begin{equation}
    \left\|
    -\frac{\boldsymbol{\xi}}{\xi_c}
    -\alpha_s\boldsymbol{s}_q
    \right\|_2
    =
    \left\|
    \frac{\boldsymbol{\xi}}{\xi_c}
    -\alpha_s\boldsymbol{s}_{q^{\ast}}
    \right\|_2.
    \label{eq:paired-source-distance}
\end{equation}
The paired points have equal weights as well. Reindexing the source sum in Eq.~\eqref{eq:source-conditioned-radial-basis} therefore gives
\begin{equation}
    G_{\ell,o,i}(-\boldsymbol{\xi};S_Q)
    =
    G_{\ell,o,i}(\boldsymbol{\xi};S_Q).
    \label{eq:source-conditioned-even-weight}
\end{equation}
The taper defined in Eq.~\eqref{eq:appendix-cosine-taper} is a function only of $\|\boldsymbol{\xi}\|_2$, making it even. With real SOG coefficients, the full source-conditioned multiplier is both real and even, as required by Eq.~\eqref{eq:weight-conjugate-symmetry}.

The centered frequency representation in the mask-to-AI model retains both $\boldsymbol{\xi}$ and $-\boldsymbol{\xi}$. Applying this real, even multiplier therefore preserves the Hermitian symmetry of real inputs and yields real outputs.

\end{document}